
\magnification=\magstep 1
\baselineskip=14pt     
\parskip=3pt plus1pt minus.5pt
\font\hd=cmbx10 scaled\magstep1
\overfullrule=0pt
\centerline{\hd Elliptic Three-folds I: Ogg-Shafarevich Theory.}
\input cyracc.def
\newfam\cyrfam
\font\tencyr=wncyr10
\def\cyr{\fam\cyrfam\tencyr\cyracc}
\def\TS{\hbox{\cyr Sh}}
\def\num{\global\advance\count10 by 1 \eqno(\the\count10)}

\def\Pthree{{\bf P^3}}

\def\Ptwo{{\bf P^2}}
\def\P{{\bf P}}

\def\O{{\cal O}}
\def\L{{\cal L}}

\def\G{{\cal G}}

\def\F{{\cal F}}

\def\E{{\cal E}}

\def\ZZ{{\cal Z}}
\def\HH{{\underline H}}
\def\Pic{\mathop{\rm Pic}}
\def\Spec{\mathop{\rm Spec}}
\def\coker{\mathop{\rm coker}}
\def\Im{\mathop{\rm Im}}
\def\Brp{\mathop{\rm Br'}}
\def\Br{\mathop{\rm Br}}
\def\char{\mathop{\rm char}}

\def\supp{\mathop{\rm Supp}}

\def\boldz{{\bf Z}}

\def\Pone{{\bf P}^1}
\def\hens{{\O_{S,\bar s}}}
\def\Gm{{\bf G}_m}
\def\Gmm#1{{\bf G}_{m,#1}}
\def\qz{{\bf Q}/{\bf Z}}

\def\exact#1#2#3{0\rightarrow#1\rightarrow#2\rightarrow#3\rightarrow0}

\def\mapright#1{\smash{
  \mathop{\longrightarrow}\limits^{#1}}}
\def\mapdown#1{\Big\downarrow
   \rlap{$\vcenter{\hbox{$\scriptstyle#1$}}$}}
\bigskip
\bigskip
\centerline{\it Igor Dolgachev\footnote{*}{This work was supported by a
grant from the National Science Foundation.} and Mark Gross\footnote{**}{This
material is based upon work supported by the North Atlantic Treaty Organization
under a Grant awarded in 1990.}}
\medskip
\centerline{\it Revised version, April, 1993.}\medskip
\centerline{Department of Mathematics}
\centerline{University of Michigan}
\centerline{Ann Arbor, Michigan 48109-1003}
\def\hens{{\O_{S,\bar s}}}
\def\Gm{{\bf G}_m}
\def\Gmm#1{{\bf G}_{m,#1}}
\def\qz{{\bf Q}/{\bf Z}}

\bigskip
\bigskip
{\hd \S 0. Introduction}

Let $S$ be a variety over an algebraically closed field $k$ of characteristic
zero, with function
field $K$, and let $A$ be an elliptic curve over $K$. The Weil-Ch\^atelet
group of $A$, $WC(A)$, is
the set of principal homogeneous spaces (torsors) of $A$ over $K$, i.e.
isomorphism classes of curves of genus 1 over $K$ which have $A$ as their
jacobian. This classifies, up to birational equivalence, elliptic fibrations
over the variety $S$ with the same jacobian. Ogg-Shafarevich theory was
developed
independently in [22] and [25] to calculate this group in the case that $S$
is a curve. In this paper, we develop certain aspects of this theory in
higher dimensions.

In \S 1, we define the Tate-Shafarevich group, which is the subset of $WC(A)$,
$\TS_S(A)$, defined as follows. Let $E\in WC(A)$ and let
$f:X\rightarrow S$ be an elliptic fibration with generic fibre equal
to $E$. Then $E\in \TS_S(A)$ if
for all $s\in S$ there exists an \'etale neighborhood
$U\rightarrow S$ of $s$ such that $X\times_S U\rightarrow U$ has a rational
section. Thus, in particular, $f$ does not have multiple fibres in codimension
one.
We show how to calculate $\TS_S(A)$ when we have a model $f:X\rightarrow S$
of some torsor $E$ of $A$ where $X$ and $S$ are smooth and
has fibres only of dimension one. We relate $\TS_S(A)$ to
the Brauer group of $X$ and $S$.

In \S 2 we consider several different questions. First, we discuss when one
can find a suitable model as above to calculate the Tate-Shafarevich group.
Secondly, we wish to refine our interpretation of the Tate-Shafarevich
group by giving conditions for when an element $E\in \TS_S(A)$ has
a model $f:X\rightarrow S$ which has no isolated multiple fibres whatsoever.
For the first question, we can construct a model for $A$ over a base which
might have to be a blowup of the original base, following [17], to obtain
a so-called Miranda model which satisfies the required hypotheses to
calculate the Tate-Shafarevich group. For the second question, it turns out
that as long as $A$ has a Miranda model over $S$,
any element $E\in \TS_S(A)$ has a model where the possible isolated
multiple fibres are severely restricted and can be understood.
We
apply Mori's theory of minimal models in this analysis. Thus we obtain
an interpretation of the Tate-Shafarevich group in
terms of multiple fibres.
However, one cannot be sure that even if $E\in\TS_S(A)$, one can find
a model $f:X\rightarrow S$ of $E$ with no isolated multiple fibres
over certain singular points of the discriminant locus.
Modulo this possibility, which only occurs at collisions of Kodaira
fibre type $IV+I_0^*$,
if $Z\subseteq S$ is a closed subset, then $\TS_{S-Z}(A)/\TS_S(A)$ should
give us the new elliptic fibrations having multiple fibres only along $Z$.

In \S 3 and \S 4, we obtain information
about this quotient in special situations.
In \S 3 we study the case that $S$ is strictly local of dimension 2 and
$Z$ is the closed point. Thus we classify possibilities for fibrations
over $S$ with isolated
multiple fibres which don't have a rational section.
We find that an elliptic threefold can indeed have
such isolated multiple fibres, but these can only occur at collision points
between two different components of the discriminant locus. In particular,
if we only allow those sorts of collisions which may occur in Miranda models,
we find that isolated multiple fibres can only occur at collisions of Kodaira
fibre type $I_{M_1}+I_{M_2}^*$, with $M_1$ even, and $I_0^*+III$. This
coincides
with local calculations performed using different techniques in [21].

In \S 4, we study a similar situation where $Z$ is a curve.
This question will be considered in more detail in the sequel to this paper.
([9]) We easily reproduce further local calculations in [21], and
also obtain some global results.
In particular, we obtain a global result of the following nature:
if $Z\subseteq S$ is a complete
irreducible curve intersecting the discriminant locus
transversally (in particular the intersection is non-empty,)
and $f^{-1}(Z)\rightarrow Z$ is an elliptic surface with
all fibres irreducible, then $\TS_{S-Z}(A)\cong \TS_S(A)$. Thus it is
impossible
to perform a `logarithmic transformation' along $Z$.

Finally, in \S 5, we apply our results to give a new example of a unirational
variety which is not rational by showing its Brauer group is
$\boldz/3\boldz$. This appears to be the first example of a unirational
three-fold with three-torsion in its Brauer group.
We construct this example by taking a general
net of cubics $S$ in $\Ptwo$ and considering the incidence variety $X$
of this net
in $S\times\Ptwo$. The morphism
$X\rightarrow S$ is then an elliptic fibration without a
section, and we can take its jacobian, $J\rightarrow S$. One then finds
that the Brauer group of $J$ is $\boldz/3\boldz$, and $J$ is unirational.
Using the same technique with a net of elliptic curves on $\Pone\times\Pone$,
one
recovers the example of Artin-Mumford [2].

We would like to thank Y. Kawamata, R. Miranda and A. Verra for their useful
remarks
concerning our work. We also thank the referee for his helpful comments.

We adopt the following notation:
{\obeylines

$S$: a normal integral excellent scheme. In \S 1, we impose no further
restrictions, but in \S 2-4, we assume $S$ is a variety of finite type over an
algebraically closed field of characteristic zero, or else an open subscheme of
a strict localization  of such a variety.
$S^{(n)}$: the set of points in $S$ of codimension $n$, i.e. such that
$\dim\O_{S,s}=n$.
$\eta$: the generic point of $S$.
$i:\eta\rightarrow S$: the natural inclusion morphism.
$K=K(S)$: the field of rational functions on $S$.
$A$: an abelian variety of dimension 1 over $K$.
$\O_{S,\bar s}$: the strict henselization of the local ring $\O_{S,s}$.
$K_{\bar s}$: the field of fractions of $\hens$.
$\eta_{\bar s}$: $\Spec K_{\bar s}$
$A(\bar s)=A\times_{\eta}\eta_{\bar s}$.
$k(s)$: the residue field of a point $s\in S$.
${}_nG$: the elements of an abelian group $G$ killed by multiplication by $n$.
$G_{tors}$: the subgroup of torsion elements of an abelian group $G$.

}

All sheaves will be in the \'etale topology, and all cohomology will be
\'etale unless otherwise stated. We use [16] as a basic reference
for \'etale cohomology and follow its notation.

{\hd \S 1. The Tate-Shafarevich Group.}

We shall identify $A$ with the sheaf in the \'etale topology of
$\eta$ which it represents. The group
$$WC(A)=H^1(\eta,A)=H^1(Gal(K_{sep}/K),A(K_{sep}))$$
is called the Weil-Ch\^atelet group of $A$, where $K_{sep}$ is the separable
closure of $K$. This is the group of
isomorphism classes of torsors (i.e. principal homogeneous spaces) of
$A$ over $K$. Each such torsor is a smooth curve $E$ of genus 1
over $K$ which has $A$ as its Jacobian $Jac(E)$.
It can be identified with its Picard scheme $Pic(E)^n$ representing
locally free sheaves of some degree $n$. The group scheme $A=Jac(E)=Pic(E)^0$
acts on $Pic(E)^n$ by translations. A torsor $E$ is trivial, i.e.
$E\cong A$, if and only if $E(K)\not=\phi$.

For any $s\in S$, closed or not, there is a natural specialization map
$$loc_{\bar s}:WC(A)\rightarrow WC(A(\bar s)),$$ taking the class of $E$ to the
class of $E\times_K K_{\bar s}$.
The {\it (geometric) Tate-Shafarevich group}, $\TS_S(A)$, is defined by setting
$$\TS_S(A)=\bigcap_{s\in S}\ker(loc_{\bar s}).$$
This consists of locally trivial torsors,
that is, all $E\in WC(A)$ such that $E(K_{\bar s})\not=\phi$ for all $s\in
S$. There is a standard cohomological interpretation of this group,
identical to that in dimension 1. (See, for example, [10] III (4.44).) It is
obtained by looking at the five-term exact sequence coming from the Leray
spectral sequence for the sheaf $A$ on $\eta$ and the inclusion morphism
$i:\eta\rightarrow S$:
$$0\rightarrow H^1(S,i_*A)\rightarrow H^1(\eta,A) \rightarrow
H^0(S,R^1i_*A) \rightarrow H^2(S,i_*A)\rightarrow H^2(\eta,A).\leqno{(1)}$$
Now, for any $s\in S$, the geometric stalk
$(R^1i_*A)_{\bar s}$ is isomorphic to $H^1(\eta_{\bar s},A(\bar s))$,
and the natural homomorphism
$$H^0(S,R^1i_*A)\rightarrow\prod_{s\in S} (R^1i_*A)_{\bar s}$$
is injective by, say, [16, II Prop. 2.10]. The composition
$H^1(\eta,A)\rightarrow H^0(S,R^1i_*A)\rightarrow (R^1i_*A)_{\bar s}$ coincides
with the
map $loc_{\bar s}$, so we obtain an exact sequence
$$0\rightarrow H^1(S,i_*A)\rightarrow WC(A)\rightarrow \prod_{s\in S}
WC(A(\bar s)).$$
Thus we have $$\TS_S(A)=H^1(S,i_*A).\leqno{(2)}$$
An element $E\in WC(A)$ maps to 0 in $WC(A(\bar s))$ if and only if
there exists an irreducible
\'etale neighborhood $U\rightarrow S$ of $s$ with field of rational
functions $K(U)$ such that
$E\times_K K(U)$ has a rational point over $K(U)$. Indeed, if the image
of $E$ in $(R^1i_*A)_{\bar s}$ is zero, that means that there is an
\'etale neighborhood $U\rightarrow S$ of $s$ such that the image
of $E$ in $H^1(K(U),E\times_K K(U))$ is zero, hence $E\times_K K(U)$ has
a rational point over $K(U)$.

One might also consider the arithmetic Tate-Shafarevich group. This is defined
as above by replacing strict henselizations with henselizations: i.e., if
$s\in S$, there is a natural map
$$loc_s:WC(A)\rightarrow WC(A(s))$$
where $A(s)=A\times_{\eta} \eta_s$, with $\eta_s$ the field of fractions
of the henselization of $\O_{S,s}$. The arithmetic Tate-Shafarevich group
is defined as
$$\TS_S^{arith}(A)=\bigcap_{s\in S-\{\eta\}}\ker(loc_s).$$
If $E\in\TS_S^{arith}(A)$, then for all $s\in S-\{\eta\}$, there is an \'etale
neighborhood $U\rightarrow S$ of $s$, with $u\in U$ mapping to $s\in S$
so that $k(s)=k(u)$, and $E\times_K K(U)$ has a rational point. This
group, however, is too small to be of interest geometrically. We wish to
classify elliptic fibrations without multiple fibres, which are those
fibrations having sections locally in the \'etale topology.
These appear in the
geometric Tate-Shafarevich group, but not necessarily in the arithmetic
Tate-Shafarevich group. We discuss the relationship between multiple fibres
and the Tate-Shafarevich group in \S 2.

Now our goal in this section
is to compute the (geometric)
Tate-Shafarevich group by using a flat projective model
$f:X\rightarrow S$ with regular $X$ and $S$, with the generic fibre
isomorphic to
$E$ for
some $E\in WC(A)$. If
$\dim S=1$, then by the theory of minimal models of surfaces, there
exists a unique, non-singular $X$ which is minimal over $S$. We
have no such canonical choice in higher dimensions, but in \S 2,
we shall discuss good choices if $\dim S=2$. We need

\proclaim Definition 1.1. Let $E\in WC(A)$. $f:X\rightarrow S$ is called
a {\it good model} for $E$ if $X_{\eta}\cong E$, and $f$ is flat
and proper with $X$ and $S$ integral and regular.

By assuming resolution of singularities (say $\char k=0$), one can always
find a model $f:X\rightarrow S$ with $X_{\eta}\cong E$, $f$ proper,
and $X$ and $S$ integral and regular. By  [11], IV, 15.4.2, if $X$ and $S$
are regular, then flatness is
equivalent to the fibres of $f$ being equidimensional. There is no reason
in general to expect the existence of such a flat model for a fixed $S$.

We now follow here the calculation of $H^1(S,i_*A)$ in [10] III \S 4,
generalizing to arbitrary dimension.

As in the theory where $\dim S=1$ we will compare $H^1(S,i_*A)$ with
information coming from the cohomological Brauer groups of $X$ and $S$.
Recall that for any scheme $Z$ the latter is defined by setting
$\Brp(Z)=H^2(Z,{\bf G}_m)$. The Brauer group $\Br(Z)$ is the group of
equivalence classes of Azumaya algebras over $Z$, and there is a canonical
injection $\Br(Z)\rightarrow \Brp(Z)_{tors}$. The two groups
coincide for example when $Z$ is the spectrum of a local Henselian ring
or is a regular scheme of dimension $\le 2$.
[10] and [16] are basic
references here.

\proclaim Proposition 1.2. If $Z$ is a noetherian scheme whose strictly
local rings are factorial (e.g. if $Z$ is regular) then $H^q(Z,\Gm)$
is torsion for $q\ge 2$.

Proof. See [10], II Prop. 1.4. $\bullet$

\proclaim Proposition 1.3. If $f:X\rightarrow S$ is a morphism
which is proper with fibres of dimension $\le 1$ with $S$ strictly
local, then $H^q(X,\Gm)_{tors}=0$ for all $q\ge 3$, and ${}_nH^2(X,\Gm)=0$
when $n$ is prime to the residue field characteristics of $X$.

Proof. This follows as in the proof of [10] III Cor. 3.2. $\bullet$

\proclaim Proposition 1.4. Let $f:X \rightarrow S$ be a flat, projective
morphism with one-dimensional fibres. Then
$$(R^if_*{\bf G}_m)_{tors}=0,\quad\forall i>1.$$
If in addition all the strictly local rings of $X$ are factorial, then
$$R^if_*{\bf G}_m=0,\quad\forall i>1.$$

Proof: Since $(R^if_*\Gm)_{\bar s}=H^i(X(\bar s),\Gm)$, where
$X(\bar s)=X\times_S \Spec \O_{S,\bar s}$, we are
reduced to Prop. 1.3, unless $i=2$ and $X$ has non-zero characteristic
residue fields. When $X$ is regular of dimension 2 the fact that
$(R^2f_*\Gm)_{tors}=0$ is a well-known theorem of M. Artin. In this
case $(R^2f_*\Gm)_{tors}=R^2f_*\Gm$. As was remarked by Grothendieck,
[10] pg. 105, to extend it to our situation we need
an extension of a theorem of Greenberg about approximation of Henselian
discrete valuation rings and the equality $\Br(X)=\Brp(X)_{tors}$. The
Artin Approximation Theorem [1] is the needed extension of Greenberg's
theorem. Now the equality $\Br(X)=\Brp(X)_{tors}$ follows from the following
result from [7]: $\Br(Z)=\Brp(Z)_{tors}$ if $Z$ is the union of two open
affine subschemes.
Since we have to show that the geometric fibres
of $R^if_*\Gm$ are equal to zero for $i>1$, we may assume that $S$ is the
spectrum of a strictly local Henselian ring. Let $\L$ be a relatively
very ample invertible sheaf on $X$. We can choose two divisors $H$ and
$H'$ from the linear system $|\L|$ with no common zeroes on the
fibre $X_0$ over the closed point of $S$. Then $H$ and $H'$ have no common
points on the whole $X$ (since the closure of a common point contains a
point on $X_0$), and thus $U=X-\supp(H)$ and $U'=X-\supp(H')$
are affine over $S$, hence affine. Applying Gabber's theorem to
$X=U\cup U'$, we obtain $\Br(X)=\Brp(X)_{tors}$. The theorem is proven.
For a self-contained proof, see [7], pg. 195. $\bullet$

\proclaim Corollary 1.5. If $f:X\rightarrow S$ is a good model, then
$f_*\Gmm{X}=\Gmm{S}$, and the Leray spectral sequence gives the exact
sequence
$$\eqalign{0\rightarrow &\Pic(S) \rightarrow \Pic(X) \rightarrow
H^0(S,R^1f_*\Gm)
\rightarrow \Brp(S)\rightarrow \Brp(X)
\rightarrow H^1(S,R^1f_*\Gm) \cr
\rightarrow &H^3(S,\Gm)
\rightarrow H^3(X,\Gm)
\rightarrow H^2(S,R^1f_*\Gm)
\rightarrow H^4(S,\Gm) \rightarrow H^4(X,\Gm).\cr}\leqno{(3)}$$

We define some important integers associated with a fibration $f:X\rightarrow
S$.

\proclaim Definition 1.6.
\item{$\delta_{\eta}$} is the positive generator of the image of the
degree map $\deg:\Pic(X_{\eta})\rightarrow \boldz$.
\item{$\delta'_{\eta}$} is the minimal positive degree of an element of
$\Pic(X_{\bar \eta})^G$ where $\bar\eta=\Spec K_{sep}$ and $G=Gal(K_{sep}/K)$.
One has $\delta'_{\eta}|\delta_{\eta}$.
\item{$\delta_S$} is the smallest degree of a divisor on $X_{\eta}$ whose
closure in $X$ is finite and flat over $S$.
\item{$\delta_s$} for $s\in S$ is the positive generator of the image of the
degree map $\deg:\Pic(X(\bar s)_{\eta_{\bar s}})\rightarrow\boldz$
where $X(\bar s)=X\times_S \Spec \O_{S,\bar s}$. Note that if $\delta_s=1$
then $E\times_{\eta} \eta_{\bar s}\cong A(\bar s)$.

\proclaim Lemma 1.7. For all $i$, any element of $\ker(H^i(S,\Gm)
\rightarrow H^i(X,\Gm))$ is killed by $\delta_S$.

Proof: Let $S'\subseteq X$ be finite and flat over $S$, with degree
over $S$ equal to $\delta_S$. We have a map $f'^*:H^i(S,\Gm)
\rightarrow H^i(S',\Gm)$, which is the composition of the homomorphisms
$$H^i(S,\Gm)\rightarrow H^i(X,\Gm) \rightarrow H^i(S',\Gm)$$
induced by the maps $f:X\rightarrow S$ and $S'\subseteq X$ respectively.

We also have the norm map, since $S'$ is finite and flat over $S$,
$$N:f'_*\Gmm{S'}\rightarrow \Gmm{S}$$
which induces a map on cohomology
$$f_*:H^i(S,f_*\Gmm{S'})\rightarrow H^i(S,\Gmm{S}).$$
By the degeneration of the Leray spectral sequence for finite maps,
we have $$H^i(S,f_*\Gmm{S'})=H^i(S',\Gmm{S'}).$$ The composition
$f'_*f'^*:H^i(S,\Gmm{S})\rightarrow H^i(S,\Gmm{S})$ is then just
multiplication by $\delta_S$, being the degree of $S'$ over $S$. This proves
the lemma. $\bullet$

\proclaim Definition 1.8.
Put $P_{X/S}=R^1f_*\Gmm{X}$.
There is a canonical morphism $P_{X/S}
\rightarrow i_*i^*P_{X/S}$, and we define $\E$ and $\F$ to be the
kernel and cokernel respectively of this morphism: i.e., there is an
exact sequence
$$0\rightarrow \E \rightarrow P_{X/S} \rightarrow i_*i^*P_{X/S}
\rightarrow \F \rightarrow 0.\leqno{(4)}$$

\proclaim Proposition 1.9. Let $f:X\rightarrow S$ be a good model,
and let $s\in S$. Then the canonical morphism on
stalks
$$(P_{X/S})_{\bar s} \rightarrow (i_*i^*P_{X/S})_{\bar s}$$
is the canonical restriction map
$$r_s:\Pic X(\bar s)\rightarrow \Pic(X_{\bar\eta_{\bar s}})^G$$
where $X(\bar s)=X\times_S \Spec \O_{S,\bar s}$ and $\bar\eta_{\bar s}
=\Spec (K_{\bar s})_{sep}$, the spectrum of the separable closure of
$K_{\bar s}$, and $G=Gal((K_{\bar s})_{sep}/K_{\bar s})$.

Proof: If $\bar s$ is a geometric point, consider the base change
$f_{\bar s}:X(\bar s)\rightarrow \Spec \hens$. Then we can compute
the stalk $(P_{X/S})_{\bar s}$ as $H^1(X(\bar s),\Gm)=\Pic X(\bar s)$.
Also, if $g:\eta_{\bar s}\rightarrow \eta $
is the canonical morphism corresponding to the extension of fields, then
$$(i_*i^*P_{X/S})_{\bar s}=H^0(\eta_{\bar s}, g^*i^*R^1f_*\Gmm{X}).$$
Now
$(g^*i^*R^1f_*\Gmm{X})_{\bar\eta_{\bar s}}$ is
a discrete $G$-module, and
$$H^0(\eta_{\bar s},g^*i^*R^1f_*\Gmm{X})=
(g^*i^*R^1f_*\Gmm{X})_{\bar \eta_{\bar s}}^G.$$
We have
$$\eqalign{(g^*i^*R^1f_*\Gmm{X})_{\bar\eta_{\bar s}}&=(R^1f_*\Gmm{X})_{\bar
\eta_{\bar s}}\cr
&=H^1(X_{\bar\eta_{\bar s}},\Gmm{X_{\bar\eta_{\bar s}}})\quad\hbox{(by [16],
III 1.17)}\cr
&=\Pic(X_{\bar\eta_{\bar s}})\cr
&=\Pic(E\times_K (K_{\bar s})_{sep}),}$$
where $E=X_{\eta}$.
Thus we have $$(i_*i^*P_{X/S})_{\bar s}=\Pic (X_{\bar\eta_{\bar s}})^G.$$
$\bullet$

\proclaim Proposition 1.10.
Let $f:X\rightarrow S$ be a good model. We have
\item{a)} $\F_{\bar s}$ is $p$-primary for any $s\in S^{(1)}$, where
$p$ is the characteristic of the residue field of $\O_{S,\bar s}$;
\item{b)} for any $s\in S^{(\ge 1)}$, $\F_{\bar s}$ is torsion, killed by
$\delta_s$.

Proof: Applying the
exact sequence (3) to $\bar f:X_{\eta_{\bar s}}\rightarrow
\eta_{\bar s}$, we obtain
$$0=\Pic(\eta_{\bar s})\rightarrow \Pic(X_{\eta_{\bar s}})
\rightarrow H^0(\eta_{\bar s}, R^1\bar f_*\Gmm{X_{\eta_{\bar s}}})
\rightarrow \Br(\eta_{\bar s}) \rightarrow \Brp(X_{\eta_{\bar s}})$$
and as in the proof of Proposition 1.9,
$H^0(\eta_{\bar s}, R^1\bar f_*\Gmm{X_{\eta_{\bar s}}})
=\Pic(X_{\bar\eta_{\bar s}})^G$, so we obtain
$$0\rightarrow \Pic(X_{\eta_{\bar s}})
\rightarrow\Pic(X_{\bar\eta_{\bar s}})^G\rightarrow\ker(\Br(\eta_{\bar s})
\rightarrow \Brp(X_{\eta_{\bar s}}))\rightarrow 0.$$
Thus, $\Pic(X_{\eta_{\bar s}})=\Pic(X_{\bar\eta_{\bar s}})^G$ if
$\Br{\eta_{\bar s}}=0$.
Nevertheless, by Lemma 1.7, $\ker(\Br(\eta_{\bar s})
\rightarrow \Brp(X_{\eta_{\bar s}}))$ is killed by $\delta_s$.

Now we can compute $\F$: the stalk
$\F_{\bar s}$ is the cokernel of $(P_{X/S})_{\bar s}
\rightarrow (i_*i^*P_{X/S})_{\bar s}$, which by Proposition 1.9
is the map
$$r_s:\Pic(X(\bar s))\rightarrow \Pic(X_{\bar\eta_{\bar s}})^G,$$
which is the restriction map of Cartier divisors from $X(\bar s)$ to
$X_{\bar\eta_{\bar s}}$. Now, since $X$ is regular,
taking
the closure of any element of $\Pic(X_{\eta_{\bar s}})$ in $X(\bar s)$ gives
a Cartier divisor, so $\Pic(X(\bar s))\rightarrow \Pic(X_{\eta_{\bar s}})$
is surjective. Thus the cokernel of $r_s$ is $\ker(\Br(\eta_{\bar
s})\rightarrow
\Brp(X_{\eta_{\bar s}}))$, which as remarked above is killed by $\delta_s$,
proving b). Part a) then follows from [10, III Cor. 1.3].
$\bullet$

\proclaim Corollary 1.11. If $f:X\rightarrow S$ is a good model, and
$\delta_{s}=1$ for all $s\in S$, then $\F=0$.

The next several propositions enable us to compute the cohomology of $\E$.

\proclaim Proposition 1.12.  Let $f:X\rightarrow S$ be a good
model. Suppose that for all $t\in S^{(1)}$
such that $X_t$ is not geometrically integral, $\overline{\{t\}}$, the
closure of $\{t\}$ in $S$, is normal. Then
$\E=\bigoplus_{t\in S^{(1)}} i_{t*}i_t^*\E,$
where $i_t:t\rightarrow S$ is the inclusion map, $t=\Spec k(t)$.
Furthermore, $i_t^*\E=0$ unless $X_t$ is not geometrically integral.

Proof:
There is a natural functorial map
$$\phi:\E\rightarrow\bigoplus_{t\in S^{(1)}} i_{t*}i_t^*\E,$$
so we need only check that this is an isomorphism on stalks.
Using the definition of $r_s$ in Proposition 1.9,
$\E_{\bar s}=\ker r_s$, which is precisely the subgroup $\Pic X(\bar s)_{fib}$
of $\Pic X(\bar s)$
generated by divisors whose support doesn't intersect the generic fibre,
or divisors
which are mapped by $f_{\bar s}$ to points of codimension one in
$\Spec \O_{S,\bar s}$. Note that for $t\in S^{(1)}$,
$i_t^*\E=0$ unless the fibre $X_t$
of $f:X\rightarrow S$ is not geometrically integral, for $(i_t^*\E)_{\bar t}
=\ker(\Pic(X(\bar t))\rightarrow \Pic(X_{\bar\eta_{\bar t}}))$ and
$\Pic(X(\bar t))_{fib}=0$ if $X_t$ is geometrically integral. Let $t_1,
\ldots, t_k$ be the points of codimension 1 in $\Spec \O_{S,\bar s}$
such that the geometric fibres
$X_{\bar t_i}$ are not geometrically integral. By the normality
assumption, there is a one-to-one correspondence between the $t_i$'s
and the elements $t\in S^{(1)}$ with $s\in \overline{\{t\}}$
and $X_t$ not geometrically
integral. Then $\E_{\bar s}
=\bigoplus_i G_i$, where $G_i$ is the subgroup of $\Pic(X(\bar s))$ spanned
by the integral divisors whose generic point projects down to $t_i$.
Now it is easy to see
that $G_i$ is isomorphic under $\phi$ to the direct summand
$(i_{t_i*}i_{t_i}^*\E)_{\bar s}$. This verifies the assertion.
$\bullet$


\proclaim Proposition 1.13. Let $f:X\rightarrow S$ be a good model, and
suppose that $\overline{\{t\}}$ is normal for all $t\in S^{(1)}$ with $X_t$
not geometrically integral. Then
\item{a)}
$$H^i(S,\E)=\bigoplus_{t\in S^{(1)}} H^i(S,i_{t*}
i_t^*\E),$$
where all but a finite number of these terms are zero.
\item{b)} Suppose $X$ and $S$ are schemes over an algebraically closed field
of characteristic zero.
Let $X^i_t$, $i=1,\ldots,n(t)$ be the irreducible components of
the fibre $X_t$, $\tilde X_t^i$ the normalization of $X^i_t$, $\tilde X^i_t
\rightarrow t_i\rightarrow t$ the Stein factorizations of
$\tilde X^i_t\rightarrow t$, $C(t)=\overline{\{t\}}$ and $C(t_i)$ the
normalization of $C(t)$ in $t_i$. Then there is an exact sequence
$$\eqalign{0\rightarrow&
H^1(S,i_{t*}i^*_t\E)\rightarrow H^1(C(t),\qz)\rightarrow
\prod_{i=1}^{n(t)} H^1(C(t_i),\qz)\cr
\rightarrow& H^2(S,i_{t*}i^*_t\E)\rightarrow H^2(C(t),\qz)\rightarrow
\prod_{i=1}^{n(t)} H^2(C(t_i),\qz)\cr
\rightarrow& H^3(S,i_{t*}i_t^*\E)\rightarrow H^3(C(t),\qz)\rightarrow
\prod_{i=1}^{n(t)} H^3(C(t_i),\qz).
\cr}$$
Here each map $H^j(C(t),\qz)\rightarrow H^j(C(t_i),\qz)$ is the composition
of the functorial homomorphism corresponding to the map $C(t_i)\rightarrow
C(t)$ and the homomorphism corresponding to the map $\qz\rightarrow\qz$ given
by $x\mapsto m_ix$, where $m_i$ is the multiplicity of $X^i_t$.

Proof: a) is clear from Proposition 1.12.

For b), we have $(i_t^*\E)_{\bar t}=\ker(\Pic X(\bar t)\rightarrow
\Pic X_{\bar\eta_{\bar t}})$ which is precisely $\coker(\boldz\rightarrow
\boldz^c)$ where $c$ is the number of components of the fibre over
$\bar t$ with the map $\boldz\rightarrow\boldz^c$ given by
$a\mapsto (\bar m_ia)$
where $\bar m_i$ is the multiplicity of the $i$th component of $X_{\bar t}$.
Now let $t'$ be a Galois extension of $t$, with $H=Gal(\bar t/t')$
acting on $(i_t^*\E)_{\bar t}$ by its action on the components of $X_{\bar t}$.
Then $(i_t^*\E)(t')=(i_t^*\E)^H_{\bar t}$, which is given by
$\coker(\boldz\rightarrow \boldz^{c(t')})$, where $c(t')$ is the number
of components of $X_t\times_t t'$. One way to describe the sheaf
$i_t^*\E$ is then as follows:
the number of irreducible components of $X^i_t\times_t t'$
is the same as the number of connected components of $\tilde X_t^i\times_t t'$
which is the same as that of $t_i \times_t t'$. Hence we can write,
if we put $\tilde t=\coprod_i t_i$, with $p_t:\tilde t\rightarrow t$ the
canonical map,
$$i_t^*\E=\coker(\phi_t:\boldz_t\rightarrow p_{t*}\boldz_{\tilde t}),$$
with the $\phi_t$ map defined as follows: $(\boldz_t)_{\bar t}=\boldz$ and
$$(p_{t*}\boldz_{\tilde t})_{\bar t}=\prod_i \boldz^{[t_i:t]}$$
by [16] II 3.5 c). We then define $\phi_t$ by sending
1 to $(m_i,\ldots,m_i)\in\boldz^{[t_i:t]}$, for each $i$, where $m_i$ is the
multiplicity of $X_t^i$. Then it is clear that $(\coker \phi_t)_{\bar t}
=\E_{\bar t}=(i_t^*\E)_{\bar t}$.

Now put
$$\bar C(t)=\coprod_{i=1}^{n(t)} C(t_i)$$
and let $q_t:\bar C(t)\rightarrow C(t)$ be the natural morphism.
On $t$, we have the sequence
$$\exact{\boldz_t}{p_{t*}\boldz_{\tilde t}}{i_t^*\E},$$
which we can push forward to obtain
$$\exact{i_{t*}\boldz_t}{i_{t*}p_{t*}\boldz_{\tilde t}}{i_{t*}i_t^*\E},$$
as $R^1i_{t*}\boldz_t=0$. Now if $j_t:C(t)\rightarrow S$ is the inclusion
map, then we can check that $i_{t*}\boldz_t=j_{t*}\boldz_{C(t)}$
and $i_{t*}p_{t*}\boldz_{\tilde t}=j_{t*}q_{t*}\boldz_{\bar C(t)}$, using
the fact that $C(t)$ and $\bar C(t)$ are normal. Thus we obtain
an exact sequence
$$\exact{j_{t*}\boldz_{C(t)}}{j_{t*}q_{t*}\boldz_{\bar C(t)}}{i_{t*}i_t^*\E}.
\leqno{(5)}$$
Now $H^i(S,j_{t*}\boldz_{C(t)})=H^i(C(t),\boldz_{C(t)})$, since
$j_{t*}$ is exact, as for any closed immersion, and
$$H^i(S,j_{t*}q_{t*}\boldz_{\bar C(t)})=H^i(\bar C(t),\boldz_{\bar C(t)}),$$
as $q_t$ is finite.

Now by the exact sequence
$$\exact{\boldz}{\bf Q}{\qz}$$
and the fact that $H^i(Z,{\bf Q}_Z)=0$ for any normal scheme $Z$, $i\ge 1$,
we have
$H^1(Z,{\bf Z})=0$ and $H^2(Z,\boldz)=H^1(Z,\qz)$, and we obtain b)
by using the long exact sequence of cohomology associated with (5).
$\bullet$

\proclaim Corollary 1.14. Let $f:X\rightarrow S$ be a good model over
an algebraically closed field of characteristic zero, and
suppose that $\overline{\{t\}}$ is normal for all $t\in S^{(1)}$ with $X_t$
not geometrically integral, and let
$$d_t=gcd\{m_i[t_i:t], i=1,\ldots,n(t)\}.$$
Then
$$\ker(H^1(C(t),\qz)\rightarrow\prod_{j=1}^{n(t)} H^1(C(t_j),\qz))$$
is killed by $d_t$. In particular, if $d_t=1$, and $C(t)$ is normal,
then we have
$$H^1(S,i_{t*}i^*_t\E)=0.$$
If furthermore $t_i=t$ for some $i$ for which $m_i=1$,
then we have exact sequences
$$\exact{H^1(C(t),\qz)}{\prod_{i=1}^{n(t)} H^1(C(t_i),\qz)}
{H^2(S,i_{t*}i^*_t\E)}$$
and
$$\exact{H^2(C(t),\qz)}{\prod_{i=1}^{n(t)}
H^2(C(t_i),\qz)} {H^3(S,i_{t*}i^*_t\E)}.$$

Proof: We have a diagram
$$\matrix{H^1(C(t),\qz)&\subseteq&H^1(t,\qz)\cr
\mapdown{}&&\mapdown{}\cr
H^1(C(t_i),\qz)&\subseteq&H^1(t_i,\qz)\cr}$$
since $C(t)$ is normal,
and the map $H^1(t,\qz)\rightarrow H^1(t_i,\qz)$ is the restriction map
multiplied by $m_i$. Now we have the trace map
$H^1(t_i,\qz)\rightarrow H^1(t,\qz)$ [16, V, 1.12] whose composition with
the above map is multiplication by $m_i[t_i:t]$. Thus
$$\ker(H^1(C(t),\qz)\rightarrow H^1(C(t_i),\qz))$$
is killed by $m_i[t_i:t]$, and the intersection of these kernels is thus killed
by $d_t$.

If $t_k=t$ for some $k$ for which $m_k=1$,
then it is clear that for any $j$,
$$\ker(H^j(C(t),\qz)\rightarrow \prod_{i=1}^{n(t)} H^j(C(t_i),\qz))=0,$$
as $H^j(C(t),\qz)\rightarrow H^j(C(t_k),\qz)$ is an isomorphism.
Hence the corollary follows. $\bullet$

\proclaim Corollary 1.15. With the same hypotheses as Corollary 1.14,
assume that $d_t=1$ for any
$t\in S^{(1)}$ and $C(t)$ is normal whenever $X_t$ is not geometrically
integral. Then
$$H^1(S,\E)=0.$$

\proclaim Proposition 1.16. Let $f:X\rightarrow S$ be a good model.
Assume $\delta_s=1$ for all $s\in S$, i.e.
$X_{\eta}\in \TS_S(A)$.
Then there are the following exact sequences:
$$\exact{
\boldz/\delta_{\eta}'\boldz}
{\TS_S(A)}{H^1(S,i_*i^*P_{X/S})},$$
$$H^1(S,\E)\rightarrow H^1(S,P_{X/S}) \rightarrow H^1(S,i_*i^*P_{X/S})
\rightarrow H^2(S,\E)\rightarrow H^2(S,P_{X/S}),$$
and
$$\Brp(S)\rightarrow \Brp(X) \rightarrow H^1(S,P_{X/S})
\rightarrow H^3(S,\Gm)\rightarrow H^3(X,\Gm)\rightarrow H^2(S,P_{X/S}).$$

Proof:
We need to compare $H^1(S,i_*A)$ with $H^1(S,i_*i^*P_{X/S})$. Now
$A=Pic(X_{\eta}/\eta)^0\subseteq Pic(X_{\eta}/\eta)$ as sheaves over $\eta$,
and $Pic(X_{\eta}/\eta)=i^*P_{X/S}$. Recall that for any finite
Galois extension $\eta'\rightarrow \eta$ with Galois group
$G=Gal(\bar\eta/\eta')$,
the value of the sheaf $Pic(X_{\eta}/\eta)$ on $\eta'$ is equal to
the group $\Pic(X_{\bar\eta})^G$. This allows us to define
the degree map of sheaves on $\eta$
$$d_{\eta}:Pic(X_{\eta}/\eta)\rightarrow \boldz_{\eta},$$
and $A=\ker d_{\eta}$. Applying $i_*$, we obtain
$$\exact{i_*A}{i_*i^*P_{X/S}}{\ZZ},\leqno{(6)}$$
where $\ZZ$ is a subsheaf of $\boldz_S=i_*\boldz_{\eta}$ such that
for any $s\in S$, one has $\ZZ_{\bar s}=\delta_s'\boldz\subseteq \boldz$,
where $\delta'_s$ is the smallest degree of an element from the group
$(i_*i^*P_{X/S})_{\bar s}=\Pic(X_{\bar\eta_{\bar s}})^G$ where
$G=Gal((K_{\bar s})_{sep}/K_{\bar s})$.
Clearly $\delta_s'|\delta_s$, and so if $\delta_s=1$ for all $s\in S$, we
must have
$\ZZ=\boldz$.

Now
$$\eqalign{&\Im(H^0(S,i_*i^*P_{X/S})\rightarrow H^0(S,\ZZ))\cr
=&\Im(H^0(S,i_*i^*P_{X/S})\rightarrow H^0(S,i_*\boldz_{\eta}))\cr
=&\Im(H^0(\eta,i^*P_{X/S})\rightarrow \boldz_{\eta}))\cr
=&\delta_\eta'\boldz,\cr}$$
by Definition 1.6. Now define
$\delta_{\eta}''$ by
$$\coker(H^0(S,\ZZ)\rightarrow H^0(S,\boldz))=\boldz/\delta_\eta''\boldz.$$
We must have $\delta_\eta''|\delta_\eta'$,  and if $\ZZ=\boldz$, then
$\delta''_{\eta}=1$. We obtain, from the
exact sequence
$$\exact{\ZZ}{\boldz}{M}$$ that $H^1(S,\ZZ)=
H^0(S,M)/(\boldz/\delta_\eta''\boldz),$ as $H^1(S,\boldz)=0$, and also
$$\coker(H^0(S,i_*i^*P_{X/S})\rightarrow H^0(S,\ZZ))=\delta_{\eta}''\boldz/
\delta_{\eta}'\boldz=
\boldz/(\delta_{\eta}'/\delta_{\eta}'')\boldz.$$
Thus we obtain an exact sequence
$$\eqalign{0\rightarrow
\boldz/(\delta_{\eta}'/\delta_{\eta}'')\boldz
\rightarrow &H^1(S,i_*A) \rightarrow H^1(S,i_*i^*P_{X/S})
\rightarrow H^0(S,M)/(\boldz/\delta_{\eta}''\boldz)\cr
\rightarrow &H^2(S,i_*A) \rightarrow H^2(S,i_*i^*P_{X/S}).\cr}$$
If $\delta_s=1$ for all $s\in S$, then $M=0$ and this reduces to the first
exact
sequence of Proposition 1.16. The second and third come from
the exact sequences (4) and (3), and Corollary 1.11. $\bullet$

\proclaim Corollary 1.17. Assume that $f:X\rightarrow S$ is a good model
with a section and
one of the following conditions is satisfied:
\item{(i)} for each $s\in S^{(1)}$ the residue field $k(s)$ is
algebraically closed;
\item{(ii)} for each $s\in S^{(1)}$ the fibre $X_s$ is geometrically
integral.
Then
$$\TS_S(A)\cong \coker(\Brp(S)\rightarrow \Brp(X)).$$

Case (i) is of course the case that $\dim S=1$, covered in [10, III \S 4].

{\it Example 1.18.} Let $S\cong\Ptwo$ be a
general base-point-free net of cubics
in $\Ptwo$ over an algebraically closed field of characteristic zero.
Then the incidence relation $X\subseteq \Ptwo\times S$
is an elliptic three-fold over $S$ via the projection $f:X\rightarrow S$,
and a $\Pone$-bundle over $\Ptwo$ via the projection $p:X\rightarrow \Ptwo$.
Since a general line $l$ in $\Ptwo$ intersects each cubic curve from
$S$ in three points, its pre-image $p^*l$ under the projection $p:X\rightarrow
\Ptwo$ is a finite cover of degree 3 over $S$.
This shows that $\delta_S=\delta_\eta=3$,
unless $\delta_{\eta}\le 2$. Let us show that the latter is impossible.
If $\delta_{\eta}\le 2$, there would be a divisor $D$ on $X$ intersecting
the general fibre $F$ of $f:X\rightarrow S$ in one or two points. But
$p:X\rightarrow \Ptwo$ being a $\Pone$-bundle, $\Pic X$ is generated by
$p^*l$ and $f^*l$, which
is a quasi-section of $p:X\rightarrow \Ptwo$. One has $p^*l.F=3$ and $f^*l.F
=0$, so any linear combination of these two divisors must intersect
$F$ in a number divisible by 3. Thus there does not exist a divisor $D$
on $X$ with $D.F=1$ or $2$.

Since $X$ is rational, $\Brp(X)=0$ ([10], III Cor. 7.3).
By Proposition 1.2, $H^3(S,\Gm)$ is torsion, and by the Kummer sequence
$$\exact{\mu_n}{\Gm}{\Gm},$$
we easily get that $H^3(S,\Gm)=0$, as $H^3(S,\mu_n)=0$ since $S\cong\Ptwo$.
Applying (3), we
infer that $H^1(S,P_{X/S})=0$. Since $S$ is general, all
fibres over points in $S^{(1)}$ are geometrically irreducible,
so $H^2(S,\E)=0$. Thus $H^1(S,i_*i^*P_{X/S})=0$. Now $\delta_{\eta}=3$,
and $\delta'_{\eta}|\delta_{\eta}$.
Since the fibres of
$X\rightarrow S$ are all reduced, for every point
$s\in S$ the localized morphism $X(\bar s)\rightarrow \Spec\O_{S,\bar s}$
has a section. This shows that $X_{\eta}$ has a rational point over
$K_{\bar s}$, i.e. represents a non-trivial element from the group
$\TS_S(A)$, and $\delta_s=1$ for all $s\in S$.
Proposition 1.16 then tells us that
$$\TS_S(A)=\boldz/\delta'_{\eta}\boldz,$$
where $A$ is the jacobian variety of $X_{\eta}$.
Thus $\delta'_{\eta}\not= 1$, so $\delta'_{\eta}=3$ and $$\TS_S(A)
=\boldz/3\boldz.$$

A curious question: what is the other non-trivial torsor from $\TS_S(A)$?

{\hd \S 2. Models for Elliptic Threefolds.}

We assume from now on that $S$ is an integral scheme of finite type over
an algebraically closed field of characteristic zero,
or an open subscheme of a strict localization
of such a scheme. Starting with Theorem 2.8, we also assume that $\dim S=2$.

\proclaim Definition 2.1. A projective morphism $f:X\rightarrow S$ is
called an {\it elliptic fibration} if its generic fibre $E$ is
a regular curve of genus one and all fibres are geometrically connected.
$f:X\rightarrow S$ is called a model for $E$.
The closed subset
$$\Delta=\{s\in S|\hbox{$X_s$ is not regular}\}$$
is called the {\it degeneration} or {\it discriminant locus}. If $t\in
S^{(1)}$,
the {\it fibre type of $t$} is the Kodaira fibre type ($_mI_n,I_n^*,II,II^*,
III,III^*,IV,IV^*$) of the central fibre of a relatively minimal model (or
N\'eron model) of $X(\bar t)=X\times_S \Spec \O_{S,\bar t}$.
A {\it collision} is a singular point of $\Delta$. The closed subset
$$\Delta^m=\{s\in S| \hbox{$f$ is not smooth at any $x\in f^{-1}(s)$}\}$$
is called the {\it multiple locus} of $f$. A fibre over a point
$s\in\Delta^m$ is called {\it multiple}. A fibre is called an
{\it isolated multiple fibre} if it is over a zero-dimensional
component of $\Delta^m$. Note that each component of a multiple fibre
must be either of dimension $>1$ or of dimension one and non-reduced
at each of its points. A {\it section} (resp. a {\it rational
section}) of $f$ is a closed subscheme $Y$ of $X$ for which the
restriction of $f$ to $Y$ is an isomorphism (resp. a birational morphism).

The first goal of this section is to determine when a good model for an
elliptic curve over $K$ exists. To construct such a good model, we begin
by studying Weierstrass fibrations and review their properties. Following
Miranda in [17], we can then resolve the singularities of Weierstrass models
under certain circumstances to obtain a good model.

First we review the definition of a Weierstrass model. See
[5,19,20] for details.
Let $\L$ a line bundle on a scheme $S$, $a\in H^0(S,\L^{\otimes 4})$,
and $b\in H^0(S,\L^{\otimes 6})$ such that $4a^3+27b^2$ is a non-zero
section of $\L^{\otimes 12}$. Let $\P=\P(\O_S\oplus \L^{\otimes-2}
\oplus \L^{\otimes-3})$, $\pi:\P\rightarrow S$ be the natural projection,
and $\O_{\P}(1)$  the tautological line bundle on $\P$.
We define the scheme $W(\L,a,b)$ as a closed subscheme of $\P$
given by the equation $Y^2Z=X^3+aXZ^2+bZ^3$
where  $X$, $Y$ and $Z$
are given by the sections of
$\O_{\bf P}(1)\otimes \L^{\otimes 2}$,
$\O_{\bf P}(1)\otimes \L^{\otimes 3}$, and
$\O_{\bf P}(1)$ which correspond to the natural injections
of $\L^{\otimes-2}$, $\L^{\otimes-3}$ and $\O_S$ into
$\pi_*\O_{\P}(1)=
\O_S\oplus\L^{\otimes -2}\oplus\L^{\otimes -3}$, respectively.

The structure morphism $W(\L,a,b)\rightarrow S$ is a flat elliptic fibration,
called a {\it Weierstrass fibration}. It has a section
$\sigma:S\rightarrow W(\L,a,b)$ defined by the $S$-point
$(X,Y,Z)=(0,1,0)$. It is easy to see that $\sigma(S)$ lies in the smooth
locus of $W(\L,a,b)$ if $S$ is regular.
We will call this section the {\it section at infinity}.
Its conormal bundle is isomorphic to $\L$.

The discriminant locus of $W(\L,a,b)\rightarrow S$ is equal to the
support of the Cartier divisor defined by the section $\Delta$
of $\L^{\otimes 12}$ given by $4a^3+27b^2$.
This gives the discriminant locus a scheme structure.
A Weierstrass fibration has
fibres which are irreducible plane cubics.
Let $C$ be the set of common zeroes
of sections $a$ and $b$. The fibres of $W(\L,a,b)\rightarrow S$ over
points of $C$ are cuspidal cubics.

Note that the construction of a Weierstrass fibration is functorial. If
$g':S'\rightarrow S$ is a morphism  of schemes, then
$$W(g^*\L,g^*(a),g^*(b))\cong W(\L,a,b)\times_S S'.$$
One can define the absolute invariant $j$ of
$W(\L,a,b)$ as a section of $\Pone_{S-C}/S-C$
given by
$j(s)=4a^3/\Delta$.

\proclaim Proposition 2.2. Let $f:W(\L,a,b)\rightarrow S$ and $f':W(L',a',b')
\rightarrow S$ be two Weierstrass models over an integral scheme $S$.
Assume that no fibres of $f$ and $f'$ are cuspidal cubics. If there is an
isomorphism between the
generic fibres of $f$ and $f'$ which preserves the section at
infinity, then there exists an isomorphism of $W(\L,a,b)$ with $W(\L',a',b')$
over $S$ which preserves the section at infinity.

Proof:
The modular invariants $j$ and $j'$ of $f$ and $f'$, respectively, coincide
at the generic point $\eta$ of $S$. Thus $j(S)$ and $j(S')$ are both
equal to the closure of the same point in the generic fibre of $\Pone_S
\rightarrow S,$ hence $j=j'$. By Proposition 5.3 of [5], the
functor of isomorphisms $Isom(f,f')$ preserving the
section is representable by a finite and non-ramified scheme $S'$ surjective
over $S$. By assumption, $S'(\eta)\not=0$. This obviously implies that
$S'=S$, and $f$ is isomorphic to $f'$ over $S$. $\bullet$

\proclaim Theorem 2.3. Let $f:X\rightarrow S$ be an elliptic fibration
which possesses a section $\sigma:S \rightarrow X$  with $X$ and $S$
regular.
Then there exists a birational $S$-morphism $g$
from $X$ to a Weierstrass fibration $f':W(\L,a,b)\rightarrow S$.
The sheaf $\L$ is
isomorphic to all of the following sheaves: $\sigma^*(\Omega^1_{X/S})$,
$f_*\omega_{X/S}$, $(R^1f_*\O_X)^{-1}$, or $\O_{\sigma(S)}(-\sigma(S))$ when
they are invertible.
The map $g$ blows down all components of fibres which do not intersect
$\sigma(S)$.

Proof: This is [20], Theorem 2.1. The theorem there is stated for complex
manifolds though in fact the proof works whenever $\char k\not=2,3$. See
also [5] and [19] for results of a similar nature. $\bullet$

\proclaim Proposition 2.4. Let $A$ be an elliptic curve over $K$ with
a rational point $\xi\in A(K)$. For any regular model $S$ of
$K$ there exists a Weierstrass fibration $W(\L,a,b)$ over $S$ with
generic fibre isomorphic to $A$. The closure of the point $\xi$
is the section at infinity of $W(\L,a,b)\rightarrow S$.

Proof: By resolution of singularities and resolution of indeterminate
points of rational maps we can find an elliptic fibration of regular
schemes $f:X\rightarrow S$ with generic fibre isomorphic to $A$.
A rational point of $A$ gives rise to a rational section of $f$. Since
$S$ is regular, it defines a section over an open subset $U$
whose complement is of codimension $\ge 2$.
Since $X$
and $S$ are regular and the restriction $f_U:X_U\rightarrow U$ of $f$
over $U$ admits a section,
we can apply Theorem 2.3 to $f_U$ to obtain a birational
$U$-morphism from $X_U$ to a Weierstrass fibration $W(\L,a,b)$ over
$U$. Since $S$ is regular, $\L$ can be extended to $S$. The sections
$a$ and $b$ also extend, since $S-U$ is codimension $\ge 2$.
We thus obtain a Weierstrass elliptic fibration
$W(\L,a,b)\rightarrow S$ birational to $f:X\rightarrow S$. The last
assertion follows from the construction of the birational map
$X\rightarrow W(\L,a,b)$. $\bullet$

\proclaim Proposition 2.5. Let $M$ be an invertible sheaf with a non-zero
section $e$. Then the Weierstrass fibrations $W(\L,a,b)$
and $W(\L\otimes M,a\otimes e^4,b\otimes e^6)$ are isomorphic over the
complement of the set of zeroes of the section $e$. The birational map
is defined by the formula:
$$(X,Y,Z)\mapsto (eX,Y,e^3Z).$$

Proof: This is a straightforward generalization of the case when $S$ is
a field. $\bullet$

\proclaim Definition 2.6. A Weierstrass fibration $W(\L,a,b)\rightarrow S$
is called {\it minimal} if there is no effective divisor $D$ such that
$div(a)\ge 4D$, $div(b)\ge 6D$.

\proclaim Proposition 2.7. Assume that $S$ is a nonsingular complex surface.
Let $f:W(\L,a,b)\rightarrow S$ be a minimal Weierstrass fibration whose
discriminant locus is a divisor with normal crossings. Then $W(\L,a,b)$
has Gorenstein rational singularities. Moreover, under these conditions,
$W(\L,a,b)$ is defined uniquely, up to isomorphism, by its generic fibre.

Proof: See [20], Corollaries 2.4 and 2.6. $\bullet$

By Proposition 2.5, every Weierstrass fibration is birationally
isomorphic to a minimal Weierstrass fibration. In particular, in the assertion
of Proposition 2.4 we may assume that $A$ is isomorphic to the
generic fibre of a minimal Weierstrass fibration. We will call this
fibration a {\it Weierstrass model} of $A$.

{}From now on in this section we assume that $S$ is an algebraic non-singular
surface over
an algebraically closed field of characteristic zero.

\proclaim Theorem 2.8. Let $W(\L,a,b)\rightarrow S$ be a Weierstrass
model. Then there exists a blowing-up $S'\rightarrow S$ with $S'$ regular, and
a
Weierstrass model $W(\L',a',b')\rightarrow S'$ birational
to $W(\L,a,b)\rightarrow S$ and a resolution of singularities $X'\rightarrow
W(\L',a',b')$ with composition $X'\rightarrow
W(\L',a',b')\rightarrow S'$ with $X'$ flat over $S'$.

Proof.  Miranda [17] has given an explicit algorithm for finding
such a resolution,
first by blowing up the base surface $S$ until the reduced
discriminant
locus $\Delta_{red}$ has simple normal crossings, and
continuing further so that only one of a small list of possible
collisions between components of $\Delta$ can occur, namely
the following possibilities:
$I_{M_1}+I_{M_2},
I_{M_1}+I_{M_2}^*,
II+IV,
II+I_0^*,
II+IV^*,
IV+I_0^*,
III+I_0^*$.
Here we use the notation for the Kodaira fibre type of fibres over the
generic points of the components that meet at a collision point.
He then
constructs a natural resolution. He points out that his
construction
is an algebraic space. The only place one leaves
the category of schemes is in using a small resolution of
an ordinary double point in resolving collisions of Kodaira
type $I_{m_1}$ and $I_{m_2}$ with $m_1$ and $m_2$ odd. However,
if one blows up such collision points on $S'$, one obtains over
the exceptional curve fibres of type $I_{m_1+m_2}$, and
since $m_1+m_2$ is even, we no longer have such a collision.
Miranda has constructed these resolutions locally,  but as he points out
in [17], Theorem 15.1, the resolutions will work
globally too. $\bullet$

\proclaim Definition 2.9. A {\it Miranda elliptic fibration} is
an elliptic fibration $f:X\rightarrow S$ such that
\item{a)} $X$ and $S$ are regular and $f$ is flat and has a section;
\item{b)} the discriminant locus $\Delta$ has simple normal crossing;
\item{c)} All collisions are of type $I_{M_1}+I_{M_2}$, $I_{M_1}+I_{M_2}^*$,
$II+IV$, $II+I_0^*$, $II+IV^*$, $IV+I_0^*$ or $III+I_0^*$.

\proclaim Corollary 2.10. Let $A$ be a one-dimensional abelian variety
over a field $K$ of transcendence degree 2 over an algebraically closed
field of characteristic zero. Then there exists a
Miranda elliptic fibration $f:X\rightarrow S$ where $S$ is some projective
model for $K$ with the generic fibre isomorphic
to $A$ (a {\it Miranda model} of $A$).

Now consider $E\in \TS_S(A)$ for some $S$ and $A$, and let $f:X\rightarrow
S$ be a model for $E$. By the definition of the Tate-Shafarevich group, for
each $s\in S$, there exists an \'etale neighborhood of $s$, $U\rightarrow S$,
such that $U\times_S X\rightarrow U$ has a rational (i.e. not everywhere
defined) section. Such a section fails to be defined only in codimension
$\ge 2$.
Now if $f:X\rightarrow S$ has a local section at a point $s\in S$,
and that section does not pass
through a singular point of $X$ on $f^{-1}(s)$, then the fibre $X_s$
is not multiple.
Thus,
in particular, if the singularities of $X$ are of codimension $\ge 3$,
we learn that $f:X\rightarrow S$ has no multiple fibres in
codimension one, but may have isolated multiple fibres. We would like
to find out when there might exist a model for $E$ without any multiple fibres.

Clearly we cannot expect any model for $E$ not
to have any isolated multiple fibres.
For example, we could blow up a one-dimensional fibre of $f$, and obtain
a fibre of dimension two, which is thus a multiple fibre.
Hence we see that we should
at least insist that our model $f:X\rightarrow S$ should be
relatively minimal in some sense. The suitable definition is Mori's:

\proclaim Definition 2.11. Let $f:X\rightarrow S$ be a projective morphism.
We say it is
{\it relatively minimal} if $X$ is ${\bf Q}$-factorial
and has only terminal singularities,
and if $C\subseteq X$ is any irreducible curve mapping to a point in $S$,
then $K_X.C\ge 0$, where $K_X$ is the canonical divisor (which is
${\bf Q}$-Cartier) of $X$ [13, 0.4.1]. See [13], \S 0 for details on terminal
singularities and canonical divisors on singular varieties. Terminal
singularities have codimension $\ge 3$ in $X$.

If $\dim X=2$, then this reduces to the old definition of relatively
minimal model: $X$ is in fact regular and no rational $-1$-curves are contained
in fibres of $f$.
Recall that the main result of [18] shows that a
relatively minimal model always exists when $\dim X=3$.
We shall study the properties of relatively minimal fibrations, and
then return to the issue of isolated multiple fibres.

\proclaim Proposition 2.12. Let $f:X\rightarrow S$ be a Miranda elliptic
fibration. Then $f$ is relatively minimal.

Proof: By [17], (15.3), $K_X$ is the pull-back of a divisor on $S$. Thus
$K_X.C=0$ for any curve $C$ contained in a fibre of $f$, and so $f$ is
relatively minimal. $\bullet$

Unlike in the case that $\dim X=2$, minimal models are not unique in
dimension 3, but are closely related by an elementary birational operation
called a flop. See [15] for definitions and details.

\proclaim Theorem 2.13.
Suppose that $f:X\rightarrow S$ and $f':X'\rightarrow S$
are relatively minimal fibrations, $\dim X,X'=3$,
and $g:X\rightarrow X'$ a birational
map over $S$. Then $g$ decomposes as a sequence of flops over $S$.

Proof: [12, Theorem 6.1] and [15, Theorem 4.9] $\bullet$

This leads to the following useful observation.

\proclaim Theorem 2.14. Let $f:X\rightarrow S$ be a relatively minimal
elliptic fibration such that $\Delta_{red}$ has no singular points.
Then any rational section is a regular section.

Proof: Let $W\rightarrow S$ be a minimal Weierstrass model of the generic
fibre $X_{\eta}$ such that its section at infinity extends our rational
section.
By [17] Theorem 7.2, the singularities of $W$ can be resolved
uniformly without changing $S$,
as the singularities are, locally, ${\bf A}^1
\times \{\hbox{Du Val}\}$,
where the Du Val singularities are the usual ones occurring in minimal
Weierstrass models for surfaces. Furthermore, such a resolution
is
a Miranda model and hence is a relatively minimal model $f':X'\rightarrow S$.
Thus there is a sequence of flops
connecting $X$ and $X'$. But because the resolution of singularities
was uniform, there are no rigid curves contained in the fibres,
and hence no flops can occur.  Thus $X$ is isomorphic to $X'$ over
$S$, hence our section is regular. $\bullet$

{\it Remark 2.15.} This shows that if $f:X\rightarrow S$ is a relatively
minimal model with a section, with $\Delta_{red}$ non-singular,
and if $X^{\#}\subseteq X$ is the maximal open subset on which $f$ is smooth,
then $(i_*A)(U)=X^{\#}(U)=Hom_S(U,X^{\#})$, if $U\rightarrow S$ is \'etale. In
addition, it shows that $f:X^{\#}\rightarrow S$ is a weak N\'eron model of
$A$. (See [3], Def. 1, pg. 74).

Before continuing,
let us relate the concept of relatively minimal and that of a minimal
Weierstrass model.

\proclaim Proposition 2.16. Let $f:X\rightarrow S$ be a relatively minimal
elliptic
fibration possessing a section,
and let $W(f_*\omega_{X/S},a,b)$
be the Weierstrass model given by Theorem 2.3, assuming $f_*\omega_{X/S}$
is invertible. Then
this is a minimal Weierstrass model.

Proof: Let $S_0\subseteq S$ be the finite set of singular points of
$\Delta_{red}$ and images of singular points of $X$,
and let $U=S-S_0$. Put $X_U=X\times_S U$, $f':X_U\rightarrow U$
and $i:U\rightarrow
S$ the inclusion. Then $i_*f'_*\omega_{X_U/U}=f_*\omega_{X/S}$, so if
we show the theorem is true for $X_U\rightarrow U$, then extending the
Weierstrass model, we obtain a minimal Weierstrass model for $f$. Thus
we can assume that $\Delta_{red}$ is smooth, and we replace $S$ by $U$
and $X$ by $X_U$.

Now take a minimal Weierstrass model $g:W=W(\L,a,b)\rightarrow S$ of the
generic
fibre $X_{\eta}$. It follows from the proof of Theorem 2.14 that $X$ can
be obtained from $W$ by minimal resolutions of singularities of local
type ${\bf A}^1\times\{\hbox{Du Val}\}$,
and we obtain
a birational morphism $\pi:X
\rightarrow W$. We have $\pi_*\omega_{X/S}=\omega_{W/S}$ since singularities
of $W$ are rational. Hence $f_*\omega_{X/S}=g_*\pi_*\omega_{X/S}
=g_*\omega_{W/S}$. By Theorem 2.3, we obtain
that the latter sheaf is isomorphic to $\L$. Now notice that two
Weierstrass fibrations $W(\L,a,b)$ and $W(\L,a',b')$ with isomorphic
generic fibres necessarily coincide. Therefore
$W(f_*\omega_{X/S},a,b)\cong W(\L,a,b)$ is minimal. $\bullet$

\proclaim Proposition 2.17.
Let $f:X\rightarrow S$ be a relatively minimal
elliptic fibration such that $X_{\eta}\in \TS_S(A)$, with $A$ being
the jacobian of $X_{\eta}$,
and let $g:W(\L,a,b)\rightarrow S$ be a minimal Weierstrass model for $A$.
Then the discriminant curves of $f$
and $g$ coincide.

Proof: Again, as in the proof of Proposition 2.16 we can remove a finite
set of points and assume that the discriminant locus
$\Delta(f)$ of $f$ is non-singular. Now, let
$s\in S$ be a point, and $U\rightarrow S$ an \'etale neighborhood of
$s$ such that $X_{\eta}\times_K K(U)$ has a rational point over $K(U)$.
Then $g: U\times_S X\rightarrow U$ is a relatively minimal fibration
with $\Delta(g)=\Delta(f)\times_S U$, and $g$ has a section,
by Theorem 2.14. Now at the same time, consider $g':
U\times_S W(\L,a,b)\rightarrow
U$. The generic fibres of $g$ and $g'$ are necessarily isomorphic, and thus
$g$ and $g'$ are birational. Since both fibrations are relatively minimal,
the fibrations themselves must be isomorphic. Thus, the discriminant loci
agree over $U$. Since $U$ was \'etale over $S$, the discriminant loci
agree over $S$, by letting $s$ vary in $S$. $\bullet$

As we have learned in this section,
we may have to change the base $S$ to construct a good
model of an elliptic curve over $K$. As the Tate-Shafarevich group depends on
the base $S$,
this may change this group. In fact, it may only
increase:

\proclaim Proposition 2.18. Let $A$ be an elliptic curve with a rational point
over a function
field $K$ of $S$ a surface, and let $S'\rightarrow S$ be the
blowing up of $S$ in a point $s$. Then there is a natural inclusion
$$\TS_S(A)\subseteq \TS_{S'}(A).$$

Proof: Both these groups sit inside $H^1(\eta,A)$. It is then clear
that if $E\in \TS_S(A)$ and $X\rightarrow S$ is a model for this, it
has locally (in the \'etale topology)
at $s$ a rational section. Pulling this back to $X\times_S S'\rightarrow S'$,
one obtains a local rational section over all points of the exceptional
curve on $S'$, and hence $E\in \TS_{S'}(A)$. $\bullet$

We now have several reasons to study isolated multiple fibres. As mentioned
earlier, it is possible that a fibration possesses a rational section,
but nevertheless has an isolated multiple fibre. In addition, in the
previous proposition, we see that a possible increase in the Tate-Shafarevich
group is due to isolated multiple fibres. Indeed, it is clear that  in the
above Proposition,
$$\TS_{S'}(A)\subseteq \TS_{S-\{s\}}(A),$$
so if $\TS_S(A)\not= \TS_{S'}(A),$ there must be an element $E\in
\TS_{S-\{s\}}(A)$ for which all models have a multiple fibre over $s$.
This motivates us to distinguish between three different types of isolated
multiple fibres:

\proclaim Definition 2.19. Let $f:X\rightarrow S$ be a relatively minimal
elliptic
fibration,  $s\in S$ a closed point, such that
$f^{-1}(s)$ is an isolated multiple fibre. Put $S(\bar s)=\Spec \O_{S,\bar s}$
and let $A$ be the jacobian of the
generic fibre of the basechange $f_{\bar s}:X(\bar s)
\rightarrow S(\bar s)$.
We call $f^{-1}(s)$
\item{1)} an {\it evanescent multiple fibre} if there exists a sequence
of blow-ups $S'\rightarrow S(\bar s)$ such that
$X(\bar s)_{\eta_{\bar s}}\in \TS_{S'}(A)-\{0\}$;
\item{2)} an {\it obstinate multiple fibre} if $X(\bar s)_{\eta_{\bar s}}$ does
not
have a point over $\eta_{\bar s}$, and it is not evanescent.
\item{3)} a {\it locally trivial multiple fibre} if $X(\bar s)_{\eta_{\bar s}}$
has
a point over $\eta_{\bar s}$.

{\it Example 2.20.} All three types of isolated multiple fibres exist. Here
we give examples of evanescent and locally trivial multiple fibres,
and in \S 3, we give examples of obstinate isolated multiple fibres.

To obtain an evanescent multiple fibre, let $S'\rightarrow S$ be the
blowing up at a point $s\in S$, with exceptional curve $E$, and let
$f':X\rightarrow S'$ be a locally trivial fibration such that
$X_E=X\times_{S'} E$ does not have a section. Then $f:X\rightarrow S$
cannot have a rational section, or else it would lift to a
rational section of $f'$ which necessarily gives a section
$X_E\rightarrow E$. In addition, it has an isolated multiple fibre
because the map $f$ is not flat at any point of $f^{-1}(s)$.

To obtain a specific example, consider the ruled surface $F_1$
along with the map $F_1\rightarrow \Pone$ giving the ruling. Let $Y$ be,
for example, an elliptic $K3$ surface over $\Pone$ without a section,
and put $X=F_1\times_{\Pone} Y$. This is an elliptic three-fold over
$F_1$, having the desired properties, taking $E$ to be the negative
section of $F_1$.

To obtain a locally trivial multiple fibre, consider over $S=\Spec k[[s,t]]$
the double cover given by
$$y^2=sx(x^3-t).$$
In $\Spec k[[s,t,x]]$, the branch locus has three components, $B_1=\{s=0\}$,
$B_2=\{x=0\}$, and $B_3=\{x^3=t\}$. The branch locus is singular along
$B_i\cap B_j$, $1\le i< j\le 3$. Blow up $B_1\cap B_3$ first, and then
the singular locus of the proper transform of the branch locus consists
of the proper transforms of $B_1\cap B_2$ and $B_2\cap B_3$, which are
now disjoint. Blow these two curves up. This provides a resolution of
singularities for the double cover. Over $s=t=0$, there are four components,
of multiplicities $2,4,6$ and $3$. The proper transform of $x=0$ in the double
cover is now a rational section, not defined at $s=t=0$, as it contains
the entire component of multiplicity $6$.

Note that $s=t=0$ is a collision point of type $IV+I_0^*$.

\proclaim Theorem 2.21. Let $j:J\rightarrow S$ be a Miranda elliptic
fibration, $S$ a non-singular surface,
and let $E\in WC(J_{\eta})$, and $f:X\rightarrow S$ a relatively
minimal fibration with $X_{\eta}=E$. Then $f:X\rightarrow S$ has no evanescent
multiple fibres. Furthermore, $f:X\rightarrow S$ can have
locally trivial and obstinate multiple fibres only over collision points.

Proof:  The statement about obstinate multiple fibre follows from Corollary
3.2.

Suppose there is an $s\in S$ such that $f^{-1}(s)$ is a locally
trivial isolated multiple fibre. Then there is an \'etale neighborhood
$U\rightarrow S$ of $s$ such that $X_U\rightarrow U$ has a rational
section $U'\subseteq X_U$. $X_U$ may not be ${\bf Q}$-factorial. If it is
not, then by [9], Cor. 4.5, there exists a birational morphism $X_U'\rightarrow
X_U$ with $X'_U$ being ${\bf Q}$-factorial with a one dimensional exceptional
locus. $X_U'\rightarrow U$ is relatively minimal, and birational to
$J_U\rightarrow U$, and hence is related to $J_U$ by a sequence of flops.

If $s\in S$ is not a collision point, then we can take $U$ to be sufficiently
small so that there are no collision points on $U$. Then as argued in the
proof of Theorem 2.14, no flops can be performed on $J_U\rightarrow U$,
and hence no small contractions can be performed
either. Thus $J_U\cong X_U$ over $U$, and
there can be no locally trivial multiple fibre over $s\in S$ when $s$ is
not a collision point.

To prove there are no evanescent fibres, replace $J\rightarrow S$
with $J(\bar s)\rightarrow \Spec \O_{S,\bar s}$, so that we can
assume that $S$ is strictly local. Thus $\TS_S(A)=0$. Now suppose
that there is an $S'\rightarrow S$ a blowing-up of $S$ such that
$\TS_{S'}(A)\not=0$. Then one can blow up $S'$ further and
assume that there is a Miranda model $J'\rightarrow S'$ of $A$.
Now we use Proposition 1.16 to calculate $\TS_{S'}(A)$. Since $\Brp$ is
a birational invariant, $\Brp(S)=\Brp(S')=0$, and $\Brp(J)=\Brp(J')=0$.
Thus $H^1(S,P_{X/S})=0$. To calculate $H^2(S',\E)$, notice that $\Delta'$,
the discriminant locus on $S'$, consists of a union of $\Pone$'s and
one or two components which are the proper transforms of components of
$\Delta$ on $S$. Since $\Delta$ was simple normal crossing to begin
with, no component of $\Delta'$ intersects more than two other components
of $\Delta'$. Thus for each $t\in S'^{(1)}$ in $\Delta'$, we have
$C(t)\cong \Pone$ and $C(t_i) \cong \Pone$, the latter because
$C(t_i)$ is ramified over $C(t)$ in at most two places. Thus $H^2(S',\E)=0$,
and so $\TS_{S'}(A)=0$. This shows that there are no evanescent fibres
over the closed point of $S$.
$\bullet$

{\it Remark 2.22.} One can perform a more detailed analysis of the collision
fibres in Miranda models to classify all possible locally trivial
multiple fibres which can occur in Theorem 2.21 over collision points. We
claim that the only such possibility is the one given in Example 2.20, at
a collision point of type $IV+I_0^*$. We analyze some of the easier cases here,
but omit the details necessary to consider all collisions.

Recall from the proof of Theorem 2.21 that if $X_U\rightarrow U$ has a locally
trivial multiple fibre at $s\in S$, then $X_U\rightarrow U$ is
obtained from a Miranda model $J_U\rightarrow U$ by sequence of flops and
then a sequence of small contractions.

First suppose that $s\in S$ is a collision point of type $I_{M_1}+I_{M_2}$.
Since $J_U$ is regular, so is $X_U'$ (the ${\bf Q}$-factorialization of
$X_U$), since flops do not add singularities ([15], Theorem 2.4). Thus, by
[21], Theorem 4.7, the fibre $(X_U')_s$ is a fibre of type $I_{M_1+M_2}$,
and all components are reduced. We certainly cannot contract every
component to obtain $X_U$, so there remains at least one reduced component.

If $s\in S$ is a collision of type $II+IV$, $II+I_0^*$  or $II+IV^*$, then by
[17], \S 12, we see that there are no rigid components of the fibre of $J_U$
over
$s$, and hence there are no flops or small contractions.

This leaves the cases $I_{M_1}+I_{M_2}^*$, $IV+I_0^*$ and $III+I_0^*$.
In these cases, flops can occur and they may change the multiplicities of
the components of the fibre. In the case $IV+I_0^*$, we can indeed eliminate
the one component of multiplicity one, obtaining the example given
in 2.20. In
the other cases, there was more than one reduced component to begin
with, and one can show that one cannot eliminate them all. We omit the details.

The end result of this analysis, if carried out, is that if $j:J\rightarrow S$
is a Miranda model with generic fibre $A$, $f:X\rightarrow S$ a relatively
minimal fibration with $X_{\eta}\in\TS_S(A)$, then $f$ can have isolated
multiple fibres at most only over collision points of type $IV+I_0^*$.
It may still be possible that one can find another minimal model for $f$,
$f':X'\rightarrow S$, which does not have any isolated multiple fibre at
all.

Summarizing 2.18 and 2.21, we obtain

\proclaim Theorem 2.23.
Let $A$ be an elliptic curve with a rational point
over a function field $K$ of dimension 2 over an algebraically closed
field $k$ with $\char k=0$ and let $j:J\rightarrow S$ be a Miranda
model for $A$, with $S$ projective.
If
$S'$ is another projective (not necessarily regular) model for $K$,
then $\TS_{S'}(A)\subseteq \TS_S(A)$.

Proof:
If $S'$ is another projective model for $K$, then we can find an $S''$
resolving the singularities of $S'$ such that $S''$ also has a birational
morphism to $S$. Then $\TS_{S'}(A)\subseteq \TS_{S''}(A)$ by Proposition
2.18, and also $\TS_S(A)\subseteq\TS_{S''}(A)$. But if $E\in
\TS_{S''}(A)$, $E\not\in \TS_S(A)$, then a model for $E$ over $S$ must have
at least one evanescent multiple fibre, contradicting Theorem 2.21. $\bullet$

Finally, we apply some results in this section to provide a calculation of
the corank of the Tate-Shafarevich group under certain hypotheses.

\proclaim Theorem 2.24. Let  $f:X\rightarrow S$ be a Miranda model, $\dim X=3$,
$\dim S=2$, $A$ the generic fibre of $f$, with $X$ and $S$ regular and
projective. Then
\item{a)} there is an exact sequence
$$\exact{{\Brp(X)\over \Brp(S)}}{\TS_S(A)}{G}$$
where $G$ is a finite group;
\item{b)} the corank $r$ of $\TS_S(A)$ is
$$b_2(X)-\rho(X)-(b_2(S)-\rho(S)),$$
where $b_2$ is the second Betti number and $\rho$ is the rank of the Picard
group, i.e., $\TS_S(A)$ is an extension of $(\qz)^r$ by a finite group.

Proof: a) Let $U\subseteq S$ be the largest open subset for which the reduced
discriminant locus is non-singular, i.e., $U$ is $S$ minus the finite number
of collision points. Then by purity for the Brauer group, $\Brp(S)=\Brp(U)$,
and $\Brp(X)=\Brp(X_U)$, where
$X_U=X\times_S U$. Let $G=\ker(H^2(S,\E)\rightarrow
H^2(S,P_{X/S}))$ and $G_U=\ker(H^2(U,\E)\rightarrow H^2(U,P_{X_U/U}))$. Then
as we have a diagram of exact sequences
$$\matrix{&&&&0&&&&\cr
&&&&\mapdown{}&&&&\cr
0&\mapright{}&{\Brp(X)\over \Brp(S)}&\mapright{}&\TS_S(A)&\mapright{}&G&
\mapright{}&0\cr
&&\mapdown{\cong}&&\mapdown{}&&\mapdown{}&&\cr
0&\mapright{}&{\Brp(X_U)\over \Brp(U)}&\mapright{}&\TS_U(A)&\mapright{}&G_U&
\mapright{}&0\cr}$$
we have $G\subseteq G_U$, and it is enough to prove
that $G_U$ is finite. Replace $X$ by $X_U$ and $S$ by $U$.

Recall that the sheaf $i_*A$ is the same as the sheaf given by $X^{\#}$.
(Remark 2.15) Let $X_0^{\#}$ be the scheme obtained from $X$ by deleting
those components of fibres which don't meet the zero section. $X_0^{\#}$
yields a subsheaf of $X^{\#}$, and we have an exact sequence
$$\exact{X_0^{\#}}{X^{\#}}{\G},$$
where $\G$ is a sheaf supported on the discriminant locus of $f$. As a sheaf
on $\Delta_{red}$, $\G$ is represented in the category of \'etale schemes
over $\Delta_{red}$ by
$$\coprod_t\coprod_{t_i} C(t_i)\rightarrow \Delta_{red},$$
in the notation of Proposition 1.13. Here the first disjoint union is over
all $t\in \Delta_{red}^{(1)}$, and the second is over the $t_i$ corresponding
to multiplicity one components of $X_t$. $\G$ is then locally constant over
$\Delta_{red}$, and thus $H^1(S,\G)$ is finite by [16], VI, Cor. 5.5.

Now let $G'$ be the image of $H^1(S,X^{\#})\rightarrow H^1(S,\G)$. We claim
that the map $\TS_S(A)=H^1(S,X^{\#})\rightarrow H^2(S,\E)$ factors
through the map $H^1(S,X^{\#})\rightarrow G'$. Indeed, if $s\in H^1(S,X^{\#})$,
$s$ can be represented as  a \v Cech cocycle on an open covering $\{U_i\}$
by sections $s_{ij}\in X^{\#}(U_{ij})$, $U_{ij}=U_i\times_S U_j$, with
$s_{ij}+s_{jk}-s_{ik}=0$, where addition is in the group law on $A$. These
sections can also be considered as divisors on $X$ whose restriction to
the generic fibre $A$ is degree zero, and then as divisors, $s_{ij}+s_{jk}-
s_{ik}$ is a divisor of $X\times_S U_{ijk}$ whose support is contained in
$f^{-1}(\Delta)$. The 2-cocycle $(s_{ij}+s_{jk}-s_{ik})_{ijk}$ is then
the image of $s$ in $H^2(S,\E)$, by construction of the boundary map in
\v Cech cohomology. Now if $s\in H^1(S,X^{\#})$ was the image of an element
in $H^1(S,X_0^{\#})$ then we could write $s$ as a cocycle $(s_{ij})$
with $s_{ij}\in X_0^{\#}(U_{ij})$. In this case $s_{ij}+s_{jk}-s_{ik}=0$
both in the group law and as a sum of divisor classes. Thus the image of
$s$ in $H^2(S,\E)$ is zero, and the map $\TS_S(A)\rightarrow H^2(S,\E)$
factors through $G'$. Thus $G\subseteq G'$, and since $G'$ is finite,
so is $G$.

b) By a), the corank of $\TS_S(A)$ is the corank of $\Brp(X)/\Brp(S)$. The
formula given then follows from [10], II Cor. 3.4 and the fact that $\Brp(S)
\rightarrow \Brp(X)$ is injective as $f$ has a section.
$\bullet$

{\hd \S 3. Isolated multiple fibres.}

We continue to work with $S$ a variety over an algebraically closed field of
characteristic zero, or an open subscheme of a strict localization of such.

Our goal in this section is to determine when obstinate multiple fibres might
occur. Suppose that $f:\bar X\rightarrow \bar S$ is a Miranda model for
$A$ an elliptic curve over $K$, and $P\in \bar S$ a closed point. Then if
$S=\bar S-\{P\}$ and $E\in\TS_S(A)$, then by Theorem 2.21, a relatively minimal
model for $E$, $f':\bar X'\rightarrow \bar S$, has no multiple fibres over
$S$, except possibly over some collision points,
and will have a multiple fibre over $P$, if $E\not\in \TS_{\bar S}(A)$.
If $E\not\in \TS_{\bar S}(A)$, then
by Theorem 2.21, the multiple fibre must be an obstinate multiple fibre.
Thus, we calculate below $\TS_S(A)$ when $\bar S$ is strictly local. In this
case $\TS_{\bar S}(A)=0$, so $\TS_S(A)$ tells us exactly where isolated
obstinate
multiple fibres may occur.

\proclaim Theorem 3.1. Suppose $\bar f:\bar X\rightarrow \bar S$
is a flat proper
morphism of regular schemes with a regular section with $\bar S$ strictly
local, of dimension $\ge 2$, with closed point $s$, and generic fibre $A$
a geometrically regular curve of genus 1. Assume furthermore that for
any $t\in\bar S^{(1)}$ such that $X_t$ is not geometrically integral,
$\overline{\{t\}}$ is normal.
Applying the notation of \S 1 to $f:X\rightarrow S$
with $S=\bar S-\{s\}$ and $X=\bar X \times_{\bar S} S$,
there is an exact sequence
$$0\rightarrow {\TS_S(A)}\rightarrow H^2(S,\E) \rightarrow
H^2(S,P_{X/S}) \rightarrow H^2(S,i_*A) \rightarrow H^3(S,\E).$$
If $\dim S>1$, then $H^2(S,P_{X/S})=0$. If $\dim S=1$, then $H^3(S,\E)=0$ and
the
groups $H^2(S,\E)$ and $H^2(S,P_{X/S})$ and the map between them
are defined by the commutative diagram of exact sequences
$$\matrix{0&\rightarrow&\bigoplus_t \qz&\rightarrow&\bigoplus_t (\qz)^{c(t)}
&\rightarrow&H^2(S,\E)&\rightarrow&0\cr
&&\mapdown{}&&\mapdown{}&&\mapdown{}&&\cr
0&\rightarrow&\qz&\rightarrow&(\qz)^c
&\rightarrow&H^2(S,P_{X/S})&\rightarrow&0\cr}$$
where the sum is over all $t\in S^{(1)}$ such that $X_t$ is not geometrically
integral, $c(t)$ is the number of components of $X_t$, and $c$ is the
number of components of the fibre $X_0$ over the closed point of $\bar S$.
The first arrow in the first row is given for a given $t$ by
$$a\mapsto (m_ir_ia),$$
where $m_i$ is the multiplicity of $X^i_t$ and $r_i$ is the ramification
of $C(t_i)\rightarrow C(t)$ at $s$. The first arrow in the second
row is given by $a\mapsto (m_ia)$ where $m_i$ is the multiplicity
of the $i$th component of $X_0$. The first vertical arrow is summation
and the second takes $X^i_t$ to the sum of components of the
central fibre of $\bar X^i_t \rightarrow C(t_i)$ taken with
multiplicities, where $\bar X^i_t$ is the normalization
of the closure of $X^i_t$ in $\bar X$.

Proof. By Proposition 1.3 applied to $\bar X\rightarrow \bar S$ and $\bar S
\rightarrow \bar S$, we see that $\Brp(\bar S)=\Brp(\bar X)=0$.
By the purity theorem for the Brauer group [10], pg. 135,
we have $\Brp(S)=0$ and $\Brp(X)=0$.
By Proposition 1.16, using the fact that $\delta_{\eta}=
\delta'_{\eta}=1$ we have $\TS_S(A)=H^1(S,i_*i^*P_{X/S})$. In addition
$H^1(S,\E)=0$ by Corollary 1.15 and $\F=0$ by Corollary 1.11. Thus
we obtain an exact sequence
$$0\rightarrow \TS_S(A) \rightarrow H^2(S,\E)
\rightarrow H^2(S,P_{X/S})\rightarrow H^2(S,i_*i^*P_{X/S})\rightarrow
H^3(S,\E)$$
from exact sequence (4) in Definition 1.8.
Also, by the sequence (6) in the proof of Proposition 1.16 and the equality
$\ZZ=\boldz$ since $\delta_{\eta}=1$, we obtain the exact sequence
$$0\rightarrow H^2(S,i_*A) \rightarrow H^2(S,i_*i^*P_{X/S}) \rightarrow
H^1(S,\qz).$$
By purity, $H^1(S,\qz)=H^1(\bar S,\qz)=0$,
so
$$H^2(S,i_*A)\cong H^2(S,i_*i^*P_{X/S}).$$
Since $X\rightarrow S$ has a section and $H^i$ is a contravariant functor,
$H^i(S,\Gm)\rightarrow H^i(X,\Gm)$ is injective for all $i$.
By Corollary 1.5
we get
$$H^2(S,P_{X/S})=H^3(X,\Gm)/H^3(S,\Gm).$$
Now $H^3(X,\Gm)$ is torsion, $X$ being regular, and so we compute
${}_nH^3(X,\Gm)$, the kernel of the multiplication by $n$ map, using the
Kummer sequence,
$${}_nH^3(X,\Gm)\cong H^3(X,\mu_n),$$
since $H^2(X,\Gm)=0$.

Now $H^i(\bar X,\mu_n)=H^i(X_0,\mu_n)$, for all $i$,
where $X_0$ is the fibre over $s$, by [16, VI, 2.7]
and $H^i(X_0,\mu_n)=0$, $i>2$, as $X_0$ is a curve over
an algebraically closed field. Thus
$$H^i(\bar X,\mu_n)=0, i>2.$$
Now using the purity theorem again, if we let $Z$ be the
finite set of singular points of the curve $Y=(X_0)_{red}$, we
have
$$\HH^j_Z(\bar X,\mu_n)=0, j\not= 2\dim \bar X$$
and hence, as earlier,
$$0=H^i(\bar X,\mu_n)=H^i(\bar X-Z,\mu_n), i=3,4.$$
Again, by purity,
$$\HH_{Y-Z}^j(\bar X-Z,\mu_n)=\cases{0,&if $j\not=2(\dim \bar X-1)=
2\dim\bar S$,\cr
(\mu_n)_{Y-Z},&if $j=2\dim\bar S$.\cr}$$
Using the exact sequence for local cohomology
$$0=H^3(\bar X-Z,\mu_n)
\rightarrow H^3(X,\mu_n)\rightarrow H^4_{Y-Z}(\bar X-Z,\mu_n)
\rightarrow H^4(\bar X-Z,\mu_n)=0,$$
we get
$$\eqalign{H^3(X,\mu_n)&\cong H^4_{Y-Z}(\bar X-Z,\mu_n)\cr
&=H^0(Y-Z, \HH^4_{Y-Z}(X-Z,\mu_n))\cr}$$
by the spectral sequence
$$H^i(Y-Z,\HH^j_{Y-Z}(\bar X-Z,\mu_n))\Rightarrow H^{i+j}_{Y-Z}(
\bar X-Z,\mu_n),$$
and that $\HH^j_{Y-Z}(\bar X-Z,\mu_n)=0$ for $j<4$. Now we have
$$H^3(X,\mu_n)\cong \cases{H^0(Y-Z,\mu_n),&if $\dim \bar S=2$,\cr
0& otherwise,\cr}$$
and $H^0(Y-Z,\mu_n)=(\mu_n)^c$, where $c$ is the number of connected
components of $Y-Z$, which is the same as the number of irreducible
components of $X_0$.
So, we get
$$H^3(X,\Gm)=H^3(X,\Gm)_{tors}=\lim_{\longrightarrow}
H^0(Y-Z,\mu_n)\cong (\qz)^c,$$
if $\dim\bar S=2$ and $0$ otherwise.

If one goes through the same calculation for $H^3(S,\Gm)$, one finds
$$H^3(S,\Gm)=\cases{\qz&if $\dim\bar S=2$,\cr
0&otherwise.\cr}$$
Thus, if $\dim\bar S>2$, we have
$$H^2(S,P_{X/S})=0;$$
if $\dim\bar S=2$, we have
$$H^2(S,P_{X/S})=\coker(\qz\rightarrow(\qz)^c)$$
where the map is given by $a\in\qz\mapsto (m_1a,\ldots,m_ca)$, where
$m_i$ is the multiplicity of the $i$th component of $Y$.

Next, we need to calculate $H^2(S,\E)$ and $H^3(S,\E)$.

By Corollary 1.14, if $t\in S^{(1)}$, $X^i_t$ the irreducible
components of $X_t$, and $X^i_t\rightarrow t_i\rightarrow t$ the Stein
factorization,
then
$$H^2(S,i_{t*}i^*_t\E)=\coker(H^1(C(t),\qz)\rightarrow \prod_{i=1}^{n(t)}
H^1(C(t_i),\qz)),$$ where $C(t)$ is the normalization of the
closure of $t$ in $S$ and
$C(t_i)$ the normalization of $C(t)$ in $t_i$. Now let
$C(t)'$ be the normalization of the closure of $C(t)$ in $\bar S$;
$C(t)'$ is a strictly local scheme, since the closure of $C(t)$ in
$\bar S$ is strictly local and geometrically unibranched. ([11] IV 18.8.6).

If $\dim S=1$, then $\dim C(t)'=1$ and $C(t)=t$. To compute $H^i(t,\qz)$,
$i=1,2$,
we just apply purity to $C(t)'$, since $C(t)'$ is regular, and
$H^i(C(t)',\qz)=0,$ $i>0$, so $H^1(t,\qz)=\qz$, $H^2(t,\qz)=0$.

Similarly, letting $C(t_i)'$ be the normalization of $C(t)'$ in $t_i$,
we use the same argument to show $H^1(t_i,\qz)=\qz$, $H^2(t_i,\qz)=0$.

Notice that if $\dim S>1$, this argument shows that $H^i(C(t),\qz)=0,$
$i=1,2$,
whenever $C(t)'$ is regular, and $H^j(C(t_i),\qz)=0$, $j=1,2$,
whenever $C(t_i)'$ is regular.
$\bullet$

\proclaim Corollary 3.2. Let $f:\bar X\rightarrow \bar S$,
$f:X\rightarrow S$ be as in Theorem 3.1, with $\dim \bar S=2$, and assume
that $f:\bar X\rightarrow \bar S$ is a Miranda model. Then $\TS_S(A)=0$
unless the closed point $s\in \bar S$ is a collision point. If
$s\in S$ is a collision point, then
$$\TS_S(A)=\cases{
0& $I_{M_1}+I_{M_2}$\cr
0& $I_{M_1}+I^*_{M_2}$, $M_1$ odd\cr
0& $II+IV$\cr
0& $II+I_0^*$\cr
0& $II+IV^*$\cr
0& $IV+I_0^*$\cr
\boldz/2\boldz& $III+I_0^*$\cr
\boldz/2\boldz& $I_{M_1}+I^*_{M_2}$, $M_1$ even\cr}$$
where $s\in \bar S$ is a collision point of the type given.

Proof: The first statement follows immediately from Theorem 3.1. For
the others,
we'll give one specific calculation: consider
the $I_2+I_0^*$ collision. The central fibre of Miranda's resolution
has six components, $f_1+f_2+2f_3+2f_4+f_5+f_6$; over the $I_2$
component there are two divisors, $C_1$ and $C_2$ whose central
fibres are $f_1+f_2+2f_3$ and $2f_4+f_5+f_6$ respectively. Over
the $I_0^*$ component there are five divisors, $D_1,\ldots,D_5$,
with central fibres $f_1,f_2,f_3+f_4,f_5$ and $f_6$ respectively.
Thus we obtain a diagram
$$\matrix{0&\rightarrow &(\qz)^2&\rightarrow& (\qz)^7& \rightarrow& H^2(S,\E)
&\rightarrow&0\cr
&&\mapdown{}&&\mapdown{}&&\mapdown{}&&\cr
0&\rightarrow&\qz&\rightarrow&(\qz)^6&\rightarrow&H^2(S,P_{X/S})&
\rightarrow&0\cr
}$$
where the first vertical arrow is summation and the second has matrix
$$\pmatrix{1&0&1&0&0&0&0\cr
1&0&0&1&0&0&0\cr
2&0&0&0&1&0&0\cr
0&2&0&0&1&0&0\cr
0&1&0&0&0&1&0\cr
0&1&0&0&0&0&1\cr}$$
The first map in the first row is given by
$$(a,b)\mapsto (a,a,b,b,2b,b,b)$$
and the first map of the second row is
$$a\mapsto (a,a,2a,2a,a,a).$$
It is then not hard to see that
$$\ker(H^2(S,\E)\rightarrow H^2(S,P_{X/S}))=\boldz/2\boldz,$$
generated by the class of
$$\pmatrix{1/2&0&1/2&1/2&0&0&0\cr}$$
in $(\qz)^7$.

The other calculations are similar, and we omit them.
$\bullet$

{\it Example 3.3.}
We can give an example of such an element of the Tate-Shafarevich group
$\TS_S(A)$ for a collision of type $I_2+I_0^*$:
over $k[[s,t]]$, consider the equation
$$y^2=s((x-1)^2-t)(x^2-t).$$
This has a singular locus consisting of two curves,
$s=(x-1)^2-t=y=0$ and $s=x^2-t=y=0$, each mapping 2-1
to the $t$-axis. Blowing these up, one obtains a smooth fibration
whose central fibre consists of a chain of three $\Pone$'s
each with multiplicity 2. This gives an elliptic fibration $X\rightarrow S=
\Spec k[[s,t]]$ with an isolated multiple fibre. Now there exists a
Miranda model $J\rightarrow S$ of the jacobian of $X\rightarrow S$. Thus,
by Theorem 2.21, the isolated multiple fibre cannot be evanescent. It
is not locally trivial: it is easy to see that $X\rightarrow S$ does
not have a rational section. Thus this is an example of an obstinate
isolated multiple fibre.


{\it Remark 3.4.} Nakayama calculates some similar groups in [21].
There, in the local situation, he removes the coordinate axes
from the base and computes the analytic analog of the Tate-Shafarevich group,
using the monodromy of the elliptic fibration.
In some cases, this gives exactly the same results as above; in
some cases this will give a bigger group, since there might
be multiple fibres away from the origin also. Thus he also obtains
the results of Theorem 4.2 below. The advantage
of his method is that it extends to any dimension, as one does
not need a model of the fibration; the disadvantage is that it
only works in the local situation.

{\it Remark 3.5.} Finally, we wish to summarize what we know about
isolated multiple fibres when $S$ is a surface.
If $j:J\rightarrow S$ is a Miranda model
for $A$, and $f:X\rightarrow S$ a relatively minimal model for
$X_{\eta}\in WC(A)$, then $f$ has no evanescent multiple fibres, and
it has two possible types of isolated multiple fibres. They may
be locally trivial, which only happens over collision points of $j$ of type
$IV+I_0^*$, by Remark 2.22, or they may be obstinate, which only happens over
collision points of $j$
of type $I_{M_1}+I_{M_2}^*$, $M_1$ even, and $III+I_0^*$, by Corollary 3.2.
Thus, if
$f:X\rightarrow S$ only has isolated multiple fibres, then $X_{\eta}\in
\TS_S(A)$ if the multiple fibres are only over collision points of type
$IV+I_0^*$. There is no difficulty in distinguishing between obstinate
and locally trivial multiple fibres.

{\hd \S 4. Obstructions to multiple fibres along curves.}

In this section, we assume that $S$ is a regular surface over an algebraically
closed field of characteristic zero, or an open subscheme of a strict
localization of such a surface.

Given any closed subset $Z\subseteq S$, we would like to calculate
the cokernel of the natural injection
$$0\rightarrow \TS_S(A)\rightarrow \TS_{S-Z}(A).$$
Recall that if $f:X\rightarrow S$ is a Miranda model for $A$, and
$E\in\TS_{S-Z}
(A)$, but $E\not\in
\TS_{S}(A)$, then a relatively minimal model for $E$ over $S$ has
multiple fibres
over some subset of $Z$.
Thus this cokernel
will classify fibrations with jacobian $A$ such that all non-locally trivial
multiple fibres occur only over points of $Z$.
In this paper, we shall give
calculations in some simple cases.

\proclaim Lemma 4.1. Let $X$ be an excellent
regular scheme over a field of characteristic zero, $Y$ be a closed
regular subscheme of codimension 1 at each point, $Y_i,i=1,\ldots,m$ its
irreducible components, $U=X-Y$. Then
there is an exact sequence
$$\eqalign{
0&\rightarrow H^2(X,\Gm)\rightarrow H^2(U,\Gm) \rightarrow
\bigoplus_{i=1}^m H^1(Y_i,\qz)\rightarrow H^3(X,\Gm)
\rightarrow H^3(U,\Gm)\cr
&\rightarrow \bigoplus_{i=1}^m H^2(Y_i,\qz)
\rightarrow H^4(X,\Gm).\cr}$$

Proof. This is [10], III, Cor.6.2, although Grothendieck forgets
to mention the hypothesis that $Y$ must be regular. One uses
the exact sequence of local cohomology
$$H^i_Y(X,\Gm)\rightarrow H^i(X,\Gm) \rightarrow H^i(U,\Gm)
\rightarrow H^{i+1}_Y(X,\Gm),$$
the ``local to global'' spectral sequence
$$H^i(Y,\HH^j_Y(\Gm))\Rightarrow H^{i+j}_Y(X,\Gm)$$
and the local calculations
$\HH^i_Y(X,\Gm)=0$ for $i\not=1$ and
$$\HH^1_Y(\Gm)=\bigoplus_{i=1}^n \boldz_{Y_i}.$$
Here, the requirement that $Y$ be non-singular is used
for $$\HH^i_Y(X,\Gm)=0, i\ge 3,$$
([10], III Theorem 6.1)
as we use the usual cohomological purity to prove this. $\bullet$

\proclaim Theorem 4.2. Let $\bar S$ be strictly local of
dimension 2, with $\bar f:\bar X\rightarrow \bar S$ a Miranda model with
a section
with discriminant locus $\Delta$ a smooth curve (or empty), and let $C$
be a curve intersecting $\Delta$ transversally at the closed point
of $\bar S$.
Then $\TS_{\bar S-C}(A)$ depends only on the fibre type
of the fibres over $\Delta$, and is given by
$$\TS_{\bar S-C}(A)=\cases{(\qz)^{\oplus 2}&$I_0$\cr
(\qz)\oplus D&$I_m$, $m\ge 1$,\cr
D&in all other cases,\cr}$$
where $D=Discr(M/rad(M))$ where $M$ is the sublattice of
$\Pic(f^{-1}(C))$ generated by the reducible fibre. Here $rad(M)$ is the
kernel of the natural map $i_M:M\rightarrow Hom_{\boldz}(M,\boldz)$
induced by the intersection pairing on $f^{-1}(C)$,
and $Disc(L)=\coker(i_L:L\rightarrow
Hom_{\boldz}(L,\boldz))$.

Proof: Put $S=\bar S-C$, $X=\bar X \times_{\bar S} S$. By Proposition 1.3,
we have $H^i(\bar X,\Gm)=H^i(\bar S,\Gm)=0$,
$i\ge 2$. Also, because $\Delta$ and $C$ meet transversally,
$f^{-1}(C)$ is smooth, and we can apply Lemma 4.1 to get
the diagram
$$\matrix{0&&0\cr
\mapdown{}&&\mapdown{}\cr
\Brp(S)&\cong&H^1(C,\qz)\cr
\mapdown{}&&\mapdown{}\cr
\Brp(X)&\cong&H^1(f^{-1}(C),\qz)\cr}$$
Since $C$ is also strictly local, $H^i(C,\qz)=0$, and
$H^i(f^{-1}(C),\qz)=H^i(X_0,\qz)$, $i=1,2$,
where $X_0$ is the fibre over the closed point
of $\bar S$.
Thus by Lemma 4.1,
$$\eqalign{&H^3(S,\Gm)=0,\cr
&\Brp(X)\cong H^1(X_0,\qz) \cong (\qz)^r,\cr}$$
with $r=2,1$ or $0$ depending on whether
$X_0$ is smooth, multiplicative (i.e. of type $I_m$ with $m\ge 1$), or additive
(i.e. any other singular fibre type). (See [4], page 290.)
Thus by Proposition 1.16,
$H^1(S,P_{X/S})=(\qz)^r$,
and if in addition all fibres of $f$ are irreducible, then
$H^2(S,\E)=0$, and we get $H^1(S,P_{X/S})
=\TS_S(A)$, in which case we are done. If the fibres are not irreducible,
we must calculate $\ker(H^2(S,\E)\rightarrow H^2(S,P_{X/S})).$

By Corollary 1.5, we have an exact sequence
$$0=H^3(S,\Gm)\rightarrow H^3(X,\Gm) \rightarrow H^2(S,P_{X/S})
\rightarrow \ker(H^4(S,\Gm)\rightarrow H^4(X,\Gm)).$$ Since
$f$ has a section, $\ker(H^4(S,\Gm)\rightarrow H^4(X,\Gm))=0$.
Now by Lemma 4.1,
$$H^2(S,P_{X/S})\cong H^3(X,\Gm) \cong H^2(X_0,\qz)\cong (\qz)^c,$$
where $c$ is the number of irreducible components of $X_0$.

Next we calculate $H^2(S,\E)$. We take $t$ to be the generic point of
$\Delta$. Then $H^2(S,\E)$ is the cokernel of the map
$$\qz\cong H^1(t,\qz)\rightarrow \bigoplus_{i=1}^c H^1(t_i,\qz)\cong (\qz)^c,$$
where $t_i$, $i=1,\ldots,n$, are the usual field extension. (In this
case, in fact, $t_i=t$.) The map is given by multiplication
by sending 1 to $(d_1m_1,\ldots,d_nm_n)$ where $d_i=[t_i:t]$,
$m_i$ is the multiplicity of the component.

The map $H^2(S,\E)\rightarrow H^2(S,P_{X/S})$ is then given by the intersection
map
$X_t^i\mapsto (X_0^j.X_t^i)$, where $X_0^j$ is the $j$th component of $X_0$.
Indeed,
an element of $H^1(\tilde t,\qz)$ corresponds to a divisor with
support on the fibre $X_0$, and we then restrict to that fibre.

It is not difficult to see now that the kernel of this map is the
desired group, $D$.  Indeed, $H^2(S,\E)=(M/rad(M))\otimes\qz$, and
if we identify $H^2(S,P_{X/S})$ with $Hom_{\boldz}(M,\boldz)\otimes\qz
=(\qz)^c$, the map $H^2(S,\E)\rightarrow H^2(S,P_{X/S})$ is then induced by the
map $$M/rad(M)\rightarrow Hom_{\boldz}(M/rad(M),\boldz)\subseteq Hom_{\boldz}
(M,\boldz)$$
induced by the intersection pairing. The claim then follows by the snake
lemma applied to the diagram of exact sequences with $L=M/rad(M)$:
$$\matrix{&&0&&&&&&\cr
&&\mapdown{}&&&&&&\cr
0&\mapright{}&L&\mapright{}&L\otimes{\bf Q}&\mapright{}&
L\otimes\qz&\mapright{}&0\cr
&&\mapdown{}&&\mapdown{\cong}&&\mapdown{}&&\cr
0&\mapright{}&Hom_{\boldz}(L,\boldz)&\mapright{}&Hom_{\boldz}(L,
\boldz)\otimes{\bf Q}&\mapright{}&
Hom_{\boldz}(L,\boldz)\otimes\qz&\mapright{}&0\cr}$$
Thus, in the case that the fibre type is additive,
we see that $\TS_{S}(A)=D$, and if it is multiplicative, we have an
exact sequence
$$\exact{\qz}{\TS_{S-C}(A)}{D}.$$ $\bullet$

{\it Remark 4.3.} The value of $D$ is well-known: (see, for example, [4],
Prop. 5.2.4 (iii) and Cor 5.2.3)
$$D=\cases{\boldz/n\boldz&$I_n$,\cr
\boldz/4\boldz&$I_{2n+1}^*$,\cr
\boldz/2\boldz\oplus\boldz/2\boldz&$I_{2n}^*$,\cr
\boldz/3\boldz&$IV, IV^*$,\cr
\boldz/2\boldz&$III, III^*$,\cr
0&$II, II^*$.\cr}$$

We now consider the global case.

\proclaim Theorem 4.4.
Let $C\subseteq S$ be a non-singular, complete
curve, not necessarily irreducible, no component of $C$ contained
in $\Delta$, and each component of $C$ intersecting $\Delta$ transversally,
and suppose $f:X\rightarrow S$ is a Miranda model. (In this case
$f^{-1}(C)$ will be smooth.)
If no component $Y$ of
$f^{-1}(C)$ is
trivial (i.e. $f|_Y:Y\rightarrow C'\subseteq C$ a trivial fibration,)
then
we have an exact sequence
$$\exact{\TS_S(A)}{\TS_{U}(A)}{K}$$
where $U=S-C$ and $K$ is a subgroup of the group
$$\bigoplus_{t\in U^{(1)}}
{\ker(H^0(q_t^{-1}(C\cap C(t)),\qz)\rightarrow H^2(\bar C(t),\qz))
\over \ker(H^0(C\cap C(t),\qz)\rightarrow H^2(C(t),\qz))},$$
where $q_t:\bar C(t)\rightarrow C(t)$ is the map from \S 1, and all other
maps are the natural ones.

Proof: As before, we calculate $\TS_{U}(A)$. By Lemma 4.1
we have a diagram
$$\matrix{0&\rightarrow&\Brp(S)&\rightarrow&\Brp(U)&
\rightarrow&H^1(C,\qz)\cr
&&\mapdown{}&&\mapdown{}&&\mapdown{f^*}\cr
0&\rightarrow&\Brp(X)&\rightarrow&\Brp(X-f^{-1}(C))&
\rightarrow&H^1(f^{-1}(C),\qz)\cr
&&\mapdown{}&&\mapdown{}&&\mapdown{}\cr
&&H^1(S,P_{X/S})&\rightarrow&H^1(U,P_{X/S})&\rightarrow&\coker(f^*)\cr}$$

We claim that $coker(f^*)$ is zero. Indeed, this follows from [4],
Cor. 5.2.2, using the fact that no component of $f^{-1}(C)$ is trivial.
So we have a surjection $H^1(S,P_{X/S})\rightarrow
H^1(U,P_{X/S})\rightarrow 0$. This gives a diagram
$$\matrix{
&&&&0&&0&&0\cr
&&&&\mapdown{}&&\mapdown{}&&\mapdown{}\cr
0&\rightarrow&H^1(S,P_{X/S})&
\rightarrow&\TS_S(A)&\rightarrow&H^2(S,\E)&\rightarrow&{H^3(X,\Gm)
\over H^3(S,\Gm)}\cr
&&\mapdown{}
&&\mapdown{}
&&\mapdown{}
&&\mapdown{}\cr
0&\rightarrow&H^1(U,P_{X/S})&\rightarrow&
\TS_{U}(A)&\rightarrow&H^2(U,\E)&\rightarrow&{H^3(X-f^{-1}(C),\Gm)
\over H^3(U,\Gm)}\cr
&&\mapdown{}
&&\mapdown{}
&&\mapdown{}
&&\cr
&&0&\rightarrow&K_1&\rightarrow&K_2&&\cr
&&&&\mapdown{}
&&\mapdown{}
&&\cr
&&&&0&&0&&\cr}$$
where $K_1$ is what we are trying to calculate. Let us calculate
$K_2$.

We have
$$K_2\cong \bigoplus_{t\in U^{(1)}} {H^2(U,i_{t*}i_t^*\E)
\over H^2(S,i_{t*}i_t^*\E)}.$$
Let $t\in U^{(1)}$. Then
$$H^2(S,i_{t*}i_t^*\E)=\coker(H^1(C(t),\qz)\rightarrow H^1(\bar C(t),\qz)),$$
and
$$H^2(U,i_{t*}i_t^*\E)=\coker(H^1(C(t)-C\cap C(t),\qz)
\rightarrow H^1(\bar C(t)-q_t^{-1}(C\cap C(t)),\qz)),$$
where $q_t:\bar C(t)\rightarrow C(t)$ is the natural projection.

Now by the exact sequence for local cohomology of the pair $(C(t),C(t)\cap C)$
and the constant sheaf $\qz$, we have an exact sequence

$$\eqalign{0\rightarrow &H^1(C(t),\qz) \rightarrow H^1(C(t)-C(t)\cap C,\qz)
\rightarrow H^0(C(t)\cap C,\qz)\cr
 \rightarrow& H^2(C(t),\qz)
\rightarrow H^2(C(t)-C(t)\cap C,\qz),\cr}$$
where the last term is zero as long as $C(t)\cap C\not= \phi$.
{}From this and the corresponding exact sequence for $\bar C(t)$, we
obtain that $K_2$ takes the form of the group given in the theorem. $\bullet$

\proclaim Corollary 4.5. Assume that $f:X\rightarrow S$ has no
reducible fibres. Then for any non-singular curve $C\subseteq S$
which intersects $\Delta$ transversally and such that each component
of $C$ intersects $\Delta$,
$$\TS_S(A)\cong \TS_{S-C}(A).$$
In particular, there is no torsor of $A$ which has multiple fibres
over the generic point of $C$ and isolated multiple fibres outside
of $C$.

Proof: This follows immediately from Theorem 4.4, as the fact that all fibres
are irreducible implies that $K$ is zero. $\bullet$

{\hd \S 5. An Application.}

Here we shall exhibit some elliptic threefolds that are counter-examples
to the L\"uroth problem in dimension 3.

{\it Example 5.1.} Consider the elliptic threefold $f:X\rightarrow S$ of
Example
1.18. We can realize the Jacobian $A$ of $X_{\eta}$ as the general fibre of a
minimal Weierstrass model $j:J\rightarrow S$. If the net chosen is
general, then the discriminant curve $\Delta$ of $f$ is a curve of
degree 12 with 24 cusps and 21 nodes. This is a classically known result, which
can be found in, for example, [6], page 7.
Over the general point of $\Delta$ the
fibre is of type $I_1$ and over the cusps it is of type $II$. By Proposition
2.17, $\Delta$ is the discriminant curve of $j$. The Weierstrass model
will only be singular over the nodes, and as these are collisions of type
$I_1+I_1$, there is no reason to expect there to be a small resolution. But if
we blow up these 21 points, we obtain a Weierstrass model over $S'$ which has a
flat regular resolution by Theorem 2.8, $j':J'\rightarrow S'$, and over the 21
exceptional curves $E_i$, we have fibres of type $I_2$. The discriminant
curve of $j'$ is the total transform $\Delta'$ of $\Delta$ under the blow-up
$S'\rightarrow S$. The fibres over the ordinary double points of $\Delta'$
are of type $I_3$. We use the notation from \S 1. The curves $C(t)$
corresponding to the irreducible components $E_i$ of $\Delta'$ are isomorphic
to $\Pone$. It is clear that $j'^{-1}(E_i)$ is the union of two irreducible
surfaces, hence $n(t)=2$ and each $C(t_i)\cong C(t),i=1,2$. Applying Corollary
1.14, we infer that $H^2(S,\E)=0$. By Proposition
1.16 applied to $j':J'\rightarrow S'$, we get
$$\Brp(J')\cong \TS_{S'}(A).$$
Now, by Proposition 2.18, $\TS_S(A)\subseteq \TS_{S'}(A)$, so $\boldz/3\boldz
\subseteq \Brp(J')$, by Example 1.18. In fact, this is an equality, for
it is easy to see that there are no evanescent multiple fibres at collisions
of type $I_1+I_1$ using Theorem 3.1. Thus
$\Brp(J')\not=0$ and hence $J'$ is not a rational variety,
as $\Brp$ is a birational invariant ([10], III 7.3) and $\Brp(\Pthree)=0$.
Now
we claim that $J$ (and hence $J'$)
is a unirational variety. To prove this it suffices to
construct a cover $T\rightarrow S$ such that $X_T=X\times_S T$ is a rational
variety and the projection $f_T:X_T\rightarrow T$ has a rational section. Then
the generic fibres of $f_T$ and of $j_T:J\times_S T\rightarrow T$ are
isomorphic elliptic curves, hence $J\times_S T$ is rational, and the
projection $J\times_S T \rightarrow J$ gives the required birational cover
of $J$.

Let us construct the cover $T\rightarrow S$. Recall that $S$ is identified
with the net of plane cubic curves. For any $s\in S$, let $F_s$ be the
corresponding cubic. Fix a line $l$ on $\Ptwo$. Let $T=\{(s,y)\in S\times l| y
\in F_s\}$ and $T\rightarrow S$ be the first projection. This is a finite
cover of degree 3. Now $X_T=\{(x,s,y)\in \Ptwo\times S \times l| x,y\in F_s\}$.
The projection $X_T\rightarrow T$, $(x,s,y)\mapsto (s,y)$, has its fibre over
a point $(s,y)\in T$ isomorphic to $F_s$, and the map $(s,y)\mapsto (y,s,y)$
is a section.
Finally the projection $X_T\rightarrow \Ptwo\times l, (x,s,y)\mapsto (x,y)$, is
a birational map, since there is only one cubic from the net passing through
a generic point $(x,y)$ of $\Ptwo\times l$. Thus $X_T$ is a rational variety
and we are done.

Thus we obtain that the jacobian fibration of a rational elliptic threefold
could be non-rational. Note that this does not happen for elliptic surfaces:
see [4], Prop. 5.6.1.

{\it Remark 5.2.} One should compare our result with the following result of
M. Van den Bergh [26]: Let $X\rightarrow |\O_{\Ptwo}(3)|$ be the universal
family for the space of all plane cubics over ${\bf C}$, and $A$
be the jacobian of its general fibre. Then its field of rational functions
over $k$ is a purely transcendental extension of ${\bf C}$.

The following example is similar to the previous one and was suggested to
us by R. Miranda.

{\it Remark 5.3.} Let $Q$ be a non-singular quadric surface in $\Pthree$ and
$S$ a general net of elliptic curves of bidegree $(2,2)$ on $Q$. Similarly to
Example 5.1, we consider the universal family $f:X\rightarrow S$. This is
an elliptic fibration with no multiple fibres, and with irreducible
discriminant
curve of degree 12 with 24 cusps and 22 nodes. An argument similar to the one
used in Example 1.18 shows that $f$ has no rational sections. It has a
2-section defined by fixing a line on $Q$. Let $A$ be the jacobian variety
of the generic fibre of $f$ and $f':J'\rightarrow S'$ be its Miranda model.
As in Example 5.1, we show that $\Brp(J')\cong \TS_{S'}(A)\supseteq \TS_S(A)
\cong \boldz/2\boldz$. The unirationality of $J'$ is proven by modifying the
argument from Example 5.1.

{\it Remark 5.4} As was pointed out to us by A. Verra, the variety $J$ which
is the jacobian of $f:X\rightarrow S$ from
the previous example is birationally isomorphic to the double solid
ramified along a quartic symmetroid, a counterexample to the L\"uroth
problem from [2]. Recall that the latter is defined as a non-singular
model of the double cover $\pi:Y\rightarrow W$, where $W$ is a general
web of quadrics in $\Pthree$ and the branch divisor is the quartic
surface of singular quadrics. If we fix a point $p\in W$ representing a
nonsingular quadric $Q$, the composition of $Y\rightarrow W$ and the
projection from $W$ to $\Ptwo$ from the point $p$ defines a rational map
whose fibres are elliptic curves. After blowing up the two pre-images of
$p$ in $Y$, we obtain an elliptic fibration $g:Y'\rightarrow \Ptwo$ with
two sections. Its fibre over a point $x$ representing a pencil $l_x$ of
quadrics $Q+\lambda Q'$ is the double cover of $l_x$ branched at the four
points represented by the four singular quadrics in the pencil. It is
well-known
that this elliptic curve is isomorphic to the base curve of the pencil
$Q\cap Q'$. Now the web $W$ cuts out on $Q$ a general net $S$ of
elliptic curves of bidegree $(2,2)$ which we may identify with the base
of $g$, so that we obtain that the non-singular fibres of $g$ are isomorphic
to the non-singular fibres of the elliptic fibration $f:X\rightarrow S$ from
Example 5.3. Now the assertion follows from:

\proclaim Lemma 5.5. Let $f:X\rightarrow S$ and $f':X'\rightarrow S$ be two
elliptic fibrations. Assume there is an open non-empty subset $U$ such that
the fibres of $f$ and $f'$ over closed points of $U$ are
isomorphic. Then the Jacobian varieties of the general fibres of $f$ and $g$
are isomorphic.

Proof: Let $K'$ be a finite extension of $K=K(S)$ such that both general
fibres $X_{\eta}$ and $X'_{\eta}$ acquire a rational point. Let $T\rightarrow
U$ be an \'etale map with $K(T)=K'$. Then the base changes $f_T:X_T\rightarrow
T$ and $f'_T:X'_T\rightarrow T$ have isomorphic fibres and both possess
rational
sections which we can make regular sections after shrinking $T$. We may also
assume that all fibres are smooth. Then we may identify the fibrations
$f_T$ and $f'_T$ with their minimal Weierstrass fibrations. Their absolute
invariant functions coincide on the set of closed points of $T$, and hence
coincide at the generic point. This shows that the generic fibres of
$f_T$ and $f'_T$ are isomorphic. Therefore $X_{\eta}\times_K K'\cong
X'_{\eta}\times_K K'$ so that $X_{\eta}$ and $X_{\eta'}$ are torsors
over the same abelian variety. $\bullet$

{\it Remark 5.6.} Using the previous Lemma we can easily verify that the
variety $J$ from Example 5.1 is birationally isomorphic to a minimal
Weierstrass fibration $W(\O_{\Ptwo}(1),a,b)\rightarrow \Ptwo,$
where $a$ and $b$ are defined as follows.
Let $S$ and $T$ be the quartic and sextic invariants of ternary forms
in 3 variables as described for example in [24].
Restricting $S$ and $T$ to the net of cubics defining our example,
we obtain homogeneous polynomials $a'$ and $b'$ in 3 variables of degree
4 and 6 respectively. Then we set $a=2{\sqrt 2}a', b=b'/3{\sqrt 3}$.
Note that the variety $J$ does not have any structure of a conic
bundle so that the standard techniques of proving non-rationality do not
apply. Indeed, conic bundles over rational surfaces have only 2-torsion in
their Brauer group ([27], Theorem 2).

{\hd Bibliography}

\item{[1]} Artin, M., ``Algebraic Approximation of Structures over
Complete Local Rings,'' {\it Publ. Math. IHES,} {\bf 36}, (1969), 23-58.
\item{[2]} Artin, M., Mumford, D., ``Some Elementary Examples of Unirational
Varieties Which Are Not Rational,'' {\it Proc. London Math. Soc.,}
{\bf 25}, (1972) 75--95.
\item{[3]} Bosch, S., L\"utkebohmert, W., Raynaud, M., {\it N\'eron Models,}
Ergebnisse der Mathematik
und ihrer Grenzgebiete, series 3, vol. 21, Springer-Verlag, (1990).
\item{[4]} Cossec, F., and Dolgachev, I., {\it Enriques Surfaces I},
Birkh\"auser, Boston, (1989).
\item{[5]} Deligne, P., ``Courbes Elliptiques: Formulaire d'apr\`es
J. Tate,'' in {\it Modular Functions in One Variable IV,} Lecture Notes
in Mathematics vol. 476, Springer-Verlag, 1975, pg. 53-74.
\item{[6]} Dolgachev, I., and Libgober, A., ``On the Fundamental Group
of the Complement to a Discriminant Variety,'' in {\it Algebraic Geometry},
Lecture notes in Mathematics, v. 862, Springer-Verlag, (1981), 1-25.
\item{[7]} Gabber, O., ``Some Theorems on Azumaya Algebras,''
in {\it Groupe de Brauer}, Lecture notes in Mathematics, v. 844,
Springer-Verlag, (1981)
pg. 129-209.
\item{[8]} Grassi, A., ``On Minimal Models of Elliptic Threefolds,''
{\it Math. Ann.} {\bf 290}, (1991) 287-301.
\item{[9]} Gross, M., ``Elliptic Threefolds II: Multiple Fibres,'' MSRI
preprint,
(1992).
\item{[10]} Grothendieck, A., ``Le Groupe de Brauer, I, II, III,''
In {\it Dix Expos\'es
sur la Cohomologie des Sch\'emas,} North-Holland, Amsterdam, 1968, 46-188.
\item{[11]} Grothendieck, A., and Dieudonn\'e, J., {\it El\'ements de
G\'eom\'etrie Alg\'ebrique}, {\it Publ. Math. IHES,} {\bf 32}, (1967).
\item{[12]} Kawamata, Y., ``Crepant Blowing-up of 3-dimensional Canonical
Singularities,'' {\it Ann. Math.} {\bf 127}, (1988), 93-163.
\item{[13]} Kawamata, Y., Matsuda, K., Matsuki, K., ``Introduction
to the Minimal Model Problem,'' {\it Proc. Sympos. Algebraic Geom.,}
Sendai, 1985. Advances in Pure Math. {\bf 10} (1987) 283-360.
\item{[14]} Kodaira, K., ``On Compact Complex Analytic Surfaces,''
I, {\it Ann. Math.} {\bf 71} (1960), 111-152;
II, {\it Ann. Math.} {\bf 77} (1963), 563-626;
III, {\it Ann. Math.} {\bf 78} (1963), 1-40.
\item{[15]} Koll\'ar, J., ``Flops,'' {\it Nagoya Math. J.}, {\bf 113},
(1989), 15-36.
\item{[16]} Milne, J. {\it \'Etale Cohomology}, Princeton Univ. Press, (1980).
\item{[17]} Miranda, R., ``Smooth Models for Elliptic Threefolds,'' in
{\it Birational Geometry of Degenerations,} Birkh\"auser, (1983) 85-133.
\item{[18]} Mori, S., ``Flip Theorem and the Existence of Minimal
Models for 3-folds,'' {\it J. Amer. Math. Soc.,} {\bf 1}, (1988), 117-253.
\item{[19]} Mumford, D., and Suominen, K., ``Introduction to the Theory of
Moduli,'' in {\it Algebraic Geometry, Oslo 1970,} Wolters-Noordhoff Press,
(1972) 171-222.
\item{[20]} Nakayama, N., ``On Weierstrass Models,'' {\it Algebraic
Geometry and Commutative Algebra in Honor of Masayoshi Nagata,} 405--431,
(1987).
\item{[21]} Nakayama, N., ``Local Structure of an Elliptic Fibration,''
Preprint, 1991, Univ. of Tokyo.
\item{[22]} Ogg, A., ``Cohomology of Abelian Varieties over Function Fields,''
{\it Ann. Math.,} {\bf 76}, 185--212, (1962).
\item{[23]} Raynaud, M., ``Caract\'eristique d'Euler-Poincar\'e d'un
Faisceau et Cohomologie des Vari\'et\'es Ab\'eliennes.'' In {\it Dix Expos\'es
sur la Cohomologie des Sch\'emas,} North-Holland, Amsterdam, 1968, 12-30.
\item{[24]} Salmon, G., {\it A Treatise on the Higher Plane Curves,}
Third edition, Chelsea Publ. Company.
\item{[25]}
Shafarevich, I., ``Principal Homogeneous Spaces over Function Fields,''
{\it AMS Translations,} {\bf 37}, 85--113, (1964).
\item{[26]} Van den Bergh, M., ``The Center of the Generic Division Algebra,''
{\it J. Algebra}, {\bf 127} (1989) 106-126.
\item{[27]} Zagorski, A., ``On Three Dimensional Conic Bundles,'' {\it
Mat. Zametki,} {\bf 21}, (1977),
745-758, (English Translation, {\it Mat. Notes},
{\bf 21}, (1977), 420-427).
\end